\documentclass[a4paper,12pt]{amsart}
\usepackage{graphicx} 
\usepackage{inputenc}
\usepackage{tikz}
\usetikzlibrary{cd}
\usetikzlibrary{arrows.meta}
\usepackage{CJK}
\usepackage[colorlinks=true,linkcolor=blue,citecolor=blue,urlcolor=blue]{hyperref}
\usepackage{geometry,url,color}
\usepackage{amsfonts,amssymb,amsmath,fancyhdr,mathrsfs}
\usepackage{amsthm,array,ragged2e,graphicx,enumitem}
\usepackage{booktabs,float,multirow,diagbox,longtable}
\usepackage{xcolor}
\usepackage{doi}
\usepackage{capt-of}
\usepackage[numbers,square]{natbib}

\allowdisplaybreaks
\justifying
\numberwithin{equation}{section}
\newtheorem{theorem}{Theorem}[section]
\newtheorem{definition}[theorem]{Definition}
\newtheorem{lemma}[theorem]{Lemma}
\newtheorem{proposition}[theorem]{Proposition}

\newtheorem{corollary}[theorem]{Corollary}

\newtheorem{example}[theorem]{Example}

\newcommand{\Aut}{\mathrm{Aut}}
\newcommand{\Hom}{\mathrm{Hom}}
\newcommand{\End}{\mathrm{End}}
\newcommand{\id}{\mathrm{id}}
\newcommand{\re}{\mathrm{Re}}
\newcommand{\Ker}{\mathrm{Ker}}
\newcommand{\im}{\mathrm{Im}}
\newcommand{\Coker}{\mathrm{Coker}}

\newcommand{\Gal}{\mathrm{Gal}}

\newcommand{\Ext}{\mathrm{Ext}}

\newcommand{\diff}{\mathrm{d}}

\newcommand{\rk}{\mathrm{rk}}

\newcommand{\Coh}{\mathrm{Coh}}

\newcommand{\Stab}{\mathrm{Stab}}

\newcommand{\linp}[1]{\left\langle#1\right\rangle}
\newcommand{\udv}{\underline{\dim}}

\newcommand{\supp}{\mathrm{supp}}
\newcommand{\Mor}{\mathrm{Mor}}

\newcommand{\add}{\mathrm{add}}

\newcommand{\Iso}{\mathrm{Iso}}
\newcommand{\cl}{\mathrm{cl}}
\newcommand{\1}{\mathbf{1}}

\newcommand{\Lie}{\mathrm{Lie}}

\newcommand{\tor}{\mathrm{tor}}
\newcommand{\ord}{\mathrm{ord}}
\newcommand{\ddiv}{\mathrm{div}}
\newcommand{\Pic}{\mathrm{Pic}}
\newcommand{\TCoh}{\mathrm{TCoh}}
\newcommand{\BCH}{\mathrm{BCH}}

\title[Joyce Transport on Elliptic Triples]{
Joyce Transport on Elliptic Triples
}
\author[Wenyu Gao]{
Wenyu Gao
}
\address{Department of Mathematical Sciences\\
Tsinghua University\\
Beijing 100084, P.~R.~China} \email{gwy23@mails.tsinghua.edu.cn}


\begin{document}
\begin{abstract}
We study a real family of stability conditions on $\mathcal{T}_C=D^b(\Mor(\Coh C))$ for an elliptic curve $C$ over a finite field. On the connected component containing this family, we determine the Lie subalgebra generated by the Joyce elements in the derived Hall algebra. At parameter zero, the coefficient of a nonzero triple in its Joyce element is equal to the derivative at zero of a polynomial that counts its flags of subobjects in its phase category. We use this method to determine all coefficients of morphisms between indecomposable sheaves in the same stable tube.
\end{abstract}
\maketitle

\section{Introduction}
In \cite{Ati1957}, Atiyah classified indecomposable vector bundles on an elliptic curve over an algebraically closed field. Over a finite field, Burban and Schiffmann \cite{BS2012} provided a parallel classification and described the Hall algebra of $\Coh C$. Their result on the equivalence between semistable sheaves of a fixed finite slope and torsion sheaves \cite[Theorem~1.1(iii)]{BS2012} allows us to treat all slopes at once.

Bradlow, Garc\'ia-Prada and Gothen \cite{BGG2004} studied how the moduli spaces of holomorphic triples of vector bundles change with a real stability parameter. In the setting of Bridgeland stability \cite{Bri2007}, Mart\'inez-Romero, Rinc\'on-Hidalgo and R\"uffer \cite{MRR2020} described the stability space of the derived category of triples over a complex curve of genus one. They reconstruct the real stability family of Bradlow, Garc\'ia-Prada and Gothen in Bridgeland's language \cite[Example 3.27]{MRR2020}. In the present paper, we prove that the family of stability conditions satisfies the support property and the finite counting condition, and they vary continuously as the real parameter changes. Mozgovoy and Schiffmann \cite{MS2020} counted semistable quiver bundles of type $A$ over finite fields by numerical class. Our
calculation gives a separate Hall coefficient for each isomorphism class of triples.

In \cite{Joy2008}, Joyce gave wall-crossing formulas for the Hall logarithms of semistable objects in abelian categories. The homogeneous components of those logarithms were called Joyce elements by Bridgeland and Toledano Laredo \cite{BT2012}. Kontsevich and Soibelman \cite[\S 2]{KS2008} expressed wall-crossing by ordered products for graded Lie algebras. Over a finite field, To\"en \cite{Toe2006} constructed the derived Hall algebra for the triangulated category of the category of perfect dg modules over a locally finite dg category. Xiao and Xu \cite{XX2008} then extended To\"en's construction to triangulated categories satisfying certain conditions. They also provided a direct proof of the associativity of the derived Hall algebras. Under the support property and the finite counting condition, Gao and Xu \cite{GX2026} work with Joyce elements indexed by a real phase and a numerical class in the uncompleted derived Hall algebra. They expressed each Joyce element at one stability condition as a finite rational Lie polynomial in those at any other in the same component \cite[Theorem 4.14]{GX2026}. Their theorem gives the constancy of the Lie subalgebra in each connected component of the spaces of numerical stability conditions. In the present paper, we determine the Lie subalgebra defined in the component containing the real parametrized stability conditions, by expressing its generators in terms of sheaves on $C$.
\subsection{Main result of the paper}
Let $C$ be an elliptic curve over $k=\mathbb{F}_q$. Write $\TCoh C=\Mor(\Coh C)$ for the category of triples $T=(\mathcal{F}\xrightarrow{f}\mathcal{G})$ of coherent sheaves, and $\mathcal{T}_C=D^b(\TCoh C)$ for its bounded derived category. We always write $\mathbb{Z}_+=\{1,2,3,\dots\}$ and $\mathbb{N}=\{0,1,2,3,\dots\}$.

For $\alpha\in\mathbb{R}$, following the $\alpha$-slope of \cite[Definition 2.1]{BGG2004}, we take the heart $\TCoh C$ and the central charge:
\[Z_\alpha(r_1,d_1,r_2,d_2)=-(d_1+d_2+\alpha r_1)+(r_1+r_2)i,\]
where $r_i,d_i$ are the rank and degree at vertex $i$. See also \cite[Example 3.37]{MRR2020}. The resulting stability conditions $\sigma_\alpha$ vary continuously with respect to $\alpha\in\mathbb{R}$ and satisfy the support property. At fixed numerical class and phase there are only finitely many semistable isomorphism classes (Theorems \ref{support}--\ref{continuity}). By \cite[Theorem 4.14]{GX2026}, the component $\Sigma_C$ containing this family has a common Lie subalgebra $\mathfrak{L}_C$ in the uncompleted derived Hall algebra $\mathcal{H}(\mathcal{T}_C)$.

At the stability condition $\sigma_0$, the nonzero phase categories in the heart are $\mathcal{B}_\mu=\Mor(\mathcal{C}_\mu)$, where $\mathcal{C}_\mu$ consists of the semistable sheaves of slope $\mu\in\mathbb{Q}\cup\{\infty\}$. At slope $\infty$ these are the torsion sheaves. Let $\Phi_A,\Phi_B,\Phi_D$ be the sums of the basis elements for $(\mathcal{F}\rightarrow0)$, $(0\rightarrow\mathcal{F})$ and $(\mathcal{F}\xrightarrow{\id}\mathcal{F})$, respectively, as $\mathcal{F}$ runs through the isomorphism classes in $\mathcal{C}_\mu$. In the completed Hall algebra of $\mathcal{B}_\mu$, the sum $\1_{\mathcal{B}_\mu}$ over all isomorphism classes in $\mathcal{B}_\mu$ has two factorizations (Theorem \ref{pentagon}):
\[\1_{\mathcal{B}_\mu}=\Phi_B*\Phi_A=\Phi_A*\Phi_D*\Phi_B.\]
The first product uses the subtriple $(0\rightarrow\mathcal{G})\subseteq T$; the second uses the kernel, image and cokernel of $f$. This is the analogue for elliptic triples of the $A_2$ pentagon identity in \cite[\S~4.6]{Bri2018}.

Let $w_\mu$ be the primitive rank‑degree vector of slope $\mu$, with $w_\infty=(0,1)$. Let $e_{\mu,m}$ be the component of $\log\1_{\mathcal{C}_\mu}$ of class $mw_\mu$. Taking the logarithm of the first product gives a finite Lie polynomial in each numerical class. It follows that $\mathfrak{L}_C$ is generated over $\mathbb{Q}$ by
\[[n]_*A_*e_{\mu,m},[n]_*B_*e_{\mu,m},\]
where $\mu\in\mathbb{Q}\cup\{\infty\},m\in\mathbb{Z}_+,n\in\mathbb{Z}$. Here $A_*$ and $B_*$ place the sheaf at the indicated vertex, and $[n]_*$ is induced by the shift by $n$ (Theorem \ref{generators of LC}).

We provide a method of flag polynomials to compute the coefficients of Joyce elements. For $0\neq T\in\mathcal{B}_\mu$, let $Z_T(t)$ be the polynomial whose value at a positive integer $m$ counts the flags $0=T_0\subseteq\cdots\subseteq T_m=T$ in $\mathcal{B}_\mu$. At the corresponding phase $\phi_\mu$, the coefficient of $u_T$ in $\epsilon^{\sigma_0}_{\phi_\mu,\cl T}$ is $Z_T'(0)$ (Proposition \ref{derivate at 0, Bmu}). It vanishes if $T$ has nonzero components in two distinct stable tubes (Theorem \ref{sigma0 formula}). We then use the method to compute the coefficients of Joyce elements of a type of triples. Suppose now that the source and target are indecomposable sheaves of lengths $a$ and $b$ in the same tube. If the image has length $l>0$, the coefficient is $\dfrac{1}{N(N+1)}$ with $N=a+b-l$, while the zero map between these sheaves has coefficient zero (Theorem \ref{epsilon Tabl}).
\subsection{Organization of the paper}
In Section 2, we describe the stable tubes of semistable sheaves on $C$ and recall the basic properties of elliptic triples.

In Section 3, we construct the family $\sigma_\alpha$ and prove the support property, the finite counting condition and continuity in $\alpha$.

In Section 4, we prove the Hall algebra pentagon at each slope and apply Joyce transport to the component $\Sigma_C$.

In Section 5, we use flag polynomials to compute the coefficients at $\sigma_0$ and identify generators of $\mathfrak{L}_C$.
\subsection{Conventions and notation}
In this paper, every subcategory is assumed to be full. For a category $\mathcal{C}$, we write $X\in\mathcal{C}$ to mean that $X$ is an object in $\mathcal{C}$, and $\mathcal{B}\subseteq\mathcal{C}$ to mean that $\mathcal{B}$ is a subcategory of $\mathcal{C}$.  We write $\Iso(\mathcal{C})$ for the set of isomorphism classes of objects in $\mathcal{C}$. For any $X\in\mathcal{C}$, we write $X\in\Iso(\mathcal{C})$ directly when there is no ambiguity. 

For an additive category $\mathcal{C}$, we use $K^b(\mathcal{C})$ to denote the bounded homotopy category of $\mathcal{C}$. For an abelian category $\mathcal{A}$, we use $D^b(\mathcal{A})$ to denote the bounded derived category of $\mathcal{A}$. For an abelian subcategory $\mathcal{B}$ of $\mathcal{A}$, we always assume that the inclusion $\mathcal{B}\hookrightarrow\mathcal{A}$ is exact.

Let $\mathcal{S}$ be a collection of objects of an additive category $\mathcal{C}$. We denote by $\add\mathcal{S}$ the additive hull of $\mathcal{S}$, namely the smallest additive subcategory of $\mathcal{C}$ containing all direct summands of finite direct sums of objects in $\mathcal{S}$. If $\mathcal{S}=\{X_1,\cdots,X_n\}$ is finite, we also write $\add(X_1,\cdots,X_n)=\add\mathcal{S}$.

For a triangulated category $\mathcal{D}$, we use $[1]$ to denote the shift functor of $\mathcal{D}$. For a triangle
\[X\rightarrow Y\rightarrow Z\rightarrow X[1],\]
we may write it as
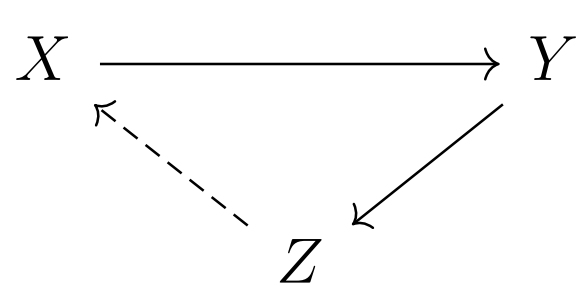
\begin{center}
    \begin{tikzcd}
    X \arrow[rr] &                      & Y \arrow[ld] \\
                 & Z \arrow[lu, dashed] &             
    \end{tikzcd}
\end{center}
when necessary. For a heart of $\mathcal{D}$, we always mean that it is a heart of a bounded $t$-structure on $\mathcal{D}$.

For an abelian or triangulated category $\mathcal{C}$, we write $K_0(\mathcal{C})$ for its Grothendieck group, whose element corresponding to an object $X\in\mathcal{C}$ is denoted by $\udv X$.
    
We always assume that $k=\mathbb{F}_q$ is a finite field with $q$ elements and fix an algebraic closure $\overline{k}$ of $k$. We write $\mathbb{P}^n$ for the $n$-dimensional projective space over $\overline{k}$. Any sheaf is assumed to be over $k$.
\section{Elliptic triples over finite fields}
 In this section, we introduce the notion of elliptic triples and recall their basic properties.
\subsection{Coherent sheaves on an elliptic curve}
We first recall the notion of elliptic curves over $k$. We only give a quick review and refer the reader to \cite{Har1977} and \cite{Sil2009} for the details. 

For an \emph{algebraic curve} $C$ over $k$, denoted by $C/k$, we mean that $C\subseteq\mathbb{P}^n$ is a $1$-dimensional projective variety over $\overline{k}$ defined by some polynomials with coefficients in $k$, where $n\in\mathbb{Z}_+$. A \emph{rational point} $x\in C$ is a point fixed by every $\sigma\in\Gal(\overline{k}/k)$, or equivalently, $x$ admits a choice of homogeneous coordinates that lie in $k$. The subset of $C$ consisting of rational points is denoted by $C(k)$. We write $|C|$ for the set of Frobenius orbits of $C$. A Frobenius orbit of $C$ is also called a \emph{closed point} on $C$. We have $C(k)\subseteq|C|$.
\begin{definition}
    An \emph{elliptic curve} over $k$ is a pair $(C,O)$, where $C$ is a nonsingular curve over $k$ with genus $1$ and $O\in C(k)$. We often omit the base point $O$ and simply call $C$ an elliptic curve over $k$.
\end{definition}

We write $\mathcal{O}_C$ for the structure sheaf over $C$. Let us recall the notion of coherent sheaves.
\begin{definition}
    A \emph{coherent sheaf} $\mathcal{F}$ over $C$ is an $\mathcal{O}_C$-module sheaf such that, for any $x\in C$, there exists a neighbourhood $U$ of $x$ together with an exact sequence of sheaves over $U$:
    \[\mathcal{O}_U^{\oplus m}\rightarrow\mathcal{O}_U^{\oplus n}\rightarrow\mathcal{F}|_U\rightarrow0,\]
    where $m,n\in\mathbb{Z}_+$ and $\mathcal{O}_U$ is the structure sheaf over $U$.
\end{definition}
We have the following examples of coherent sheaves.
\begin{example}
    \begin{enumerate}
        \item A \emph{vector bundle} over $C$ is a sheaf $\mathcal{F}$ over $C$ such that, for any $x\in C$, there exists a neighbourhood $U$ of $x$, such that
        \[\mathcal{F}|_U\cong\mathcal{O}_U^{\oplus r},\]
        where $r\in\mathbb{N}$ is independent of the choice of $x$. We call $r$ the \emph{rank} of $\mathcal{F}$ and denote it by $\rk\mathcal{F}$. It is clear that vector bundles are coherent.
        \item A \emph{torsion sheaf} over $C$ is a coherent sheaf $\mathcal{F}$ over $C$ such that its support
        \[\supp\mathcal{F}=\{x\in |C||\mathcal{F}_x\neq0\}\]
        is finite. For $x\in |C|$, we write $\mathcal{I}_x$ for the \emph{ideal sheaf} of $\mathcal{O}_C$ at $x$, that is
        \[\mathcal{I}_x(U)=\begin{cases}
            \{f\in\mathcal{O}_C(U)|f_x\in\mathfrak{m}_x\}&,x\in |U|,\\
            \mathcal{O}_{C}(U)&,x\notin |U|,
        \end{cases}\]
        where $\mathfrak{m}_x$ is the maximal ideal of the local ring $\mathcal{O}_{C,x}$. For $n\in\mathbb{Z}_+$, the quotient sheaf $\mathcal{O}_C/\mathcal{I}_x^n$ is coherent by definition. When $n=1$, we call the quotient sheaf $\underline{k}_x=\mathcal{O}_C/\mathcal{I}_x$ the \emph{skyscraper sheaf} at $x$.
        By \cite[\S1.1;\S4.1--4.2]{BS2012}, any torsion sheaf $\mathcal{F}$ has a decomposition
        \[\mathcal{F}\cong\bigoplus_{x\in\supp\mathcal{F}}\bigoplus_{j=1}^{r_x}\mathcal{O}_C/\mathcal{I}_x^{n_{x,j}},n_{x,j}\in\mathbb{Z}_+.\]
        For $x\in\supp\mathcal{F}$, we call
        \[\ell(\mathcal{F},x)=\sum_{j=1}^{r_x}n_{x,j}\]
        the \emph{length} of $\mathcal{F}$ at $x$.
    \end{enumerate}
\end{example}
We write $\Coh C$ for the category of coherent sheaves over $C$, whose morphisms are the morphisms of $\mathcal{O}_C$-modules. Then we have
\begin{theorem}[\cite{StacksProject}]
    $\Coh C$ is a hereditary abelian category over $k$.
\end{theorem}
For any $\mathcal{F}\in\Coh C$, let $\mathcal{F}_\tor$ be the maximal torsion subsheaf of $\mathcal{F}$; that is, for any open subset $U\subseteq C$, 
\[\mathcal{F}_\tor(U)=\{s\in\mathcal{F}(U)|s_x\text{ is torsion in }\mathcal{F}_x\text{ for every }x\in U\}.\]
When $\mathcal{F}$ is a torsion sheaf, we have $\mathcal{F}=\mathcal{F}_\tor$.
We have a short exact sequence
\[0\rightarrow\mathcal{F}_\tor\rightarrow\mathcal{F}\rightarrow\mathcal{F}/\mathcal{F}_\tor\rightarrow0\]
in $\Coh C$. Since $C$ is nonsingular, the quotient sheaf $\mathcal{F}/\mathcal{F}_\tor$ must be a vector bundle, by \cite[Chapter I, Theorem 6.2A; Chapter II, Ex. 5.7(b)]{Har1977}. Thus we can define rank for any coherent sheaf:
\[\rk\mathcal{F}:=\rk(\mathcal{F}/\mathcal{F}_\tor).\]
One can easily show that $\rk\mathcal{F}=0$ if and only if $\mathcal{F}$ is torsion.
\subsection{Line bundles and divisors}
Let $(C,O)$ be an elliptic curve over $k$. A \emph{line bundle} over $C$ is a vector bundle over $C$ of rank $1$. An effective tool to study line bundles over elliptic curves is the divisor.
\begin{definition}
    A \emph{divisor} on $C$ is a finite integral linear formal combination
    \begin{equation}\label{divisor}
        D=\sum_{x\in |C|}n_x x,
    \end{equation}
    where $n_x\in\mathbb{Z}$ is zero for all but finitely many $x\in |C|$.
\end{definition}
A divisor (\ref{divisor}) is called \emph{effective} if $n_x\geqslant0$ for any $x\in |C|$. Let $f\in k(C)^\times$ be a nonzero rational function. We write $\ord_xf$ to mean the order of $f$ at $x$, see \cite[\S II.1]{Sil2009}. The \emph{principal divisor} of $f$ is defined to be
\[\ddiv f=\sum_{x\in |C|}(\ord_xf)x.\]
Two divisors $D$ and $D'$ on $C$ are said to be \emph{linearly equivalent} if $D-D'=\ddiv f$ for some rational function $f$.

Given a divisor (\ref{divisor}), we construct a sheaf $\mathcal{O}_C(D)$ over $C$. For any nonempty open subset $U\subseteq C$, we define
\[\mathcal{O}_C(D)(U)=\{f\in k(C)^\times|\ord_x(f)+n_x\geqslant0,x\in |U|\}\cup\{0\}.\]
It is easy to see that if $D'$ is an effective divisor on $C$, then for any divisor $D$ on $C$, we have
\[\mathcal{O}_C(D)\subseteq\mathcal{O}_C(D+D').\]
\begin{lemma}\label{line bundles and divisors}
    $\mathcal{O}_C(D)$ is a line bundle over $C$. Moreover,
    \begin{enumerate}
        \item every line bundle over $C$ is isomorphic to $\mathcal{O}_C(D)$ for some divisor $D$ on $C$;
        \item two line bundles $\mathcal{O}_C(D)$ and $\mathcal{O}_C(D')$ are isomorphic if and only if $D$ and $D'$ are linearly equivalent.
    \end{enumerate}
\end{lemma}
\begin{proof}
    Since $C$ is a smooth curve, the result follows from \cite[Chapter II, Proposition 6.13 and Corollary 6.16]{Har1977}.
\end{proof}
Now we can define the \emph{degree} of a line bundle:
\begin{equation}\label{line degree}
    \deg\mathcal{O}_C(D)=\deg D=\sum_{x\in |C|}n_x[k_x:k],
\end{equation}
where $k_x=\mathcal{O}_{C,x}/\mathfrak{m}_x$ is the residue field at $x$.
\subsection{Picard groups and Grothendieck groups}
Let $(C,O)$ be an elliptic curve over $k$. We can use divisors to characterize the tensor products and duals of line bundles.
\begin{lemma}[{\cite[Chapter II, Propositions 6.12 and 6.13(b)]{Har1977}}]
    Let $D,D'$ be two divisors on $C$. Then
    \[\mathcal{O}_C(D)\otimes\mathcal{O}_C(D')\cong\mathcal{O}_C(D+D'),\mathcal{O}_C(D)^\vee\cong\mathcal{O}_C(-D).\]
\end{lemma}
As a consequence, if we set $\Pic(C)$ to be the set of isomorphism classes of line bundles over $C$, then the tensor product induces a group structure on $\Pic(C)$, where the inverses are induced by duals. We call the resulting group $(\Pic(C),\otimes,^\vee)$ the \emph{Picard group} of $C$. It follows from Lemma \ref{line bundles and divisors} that $\Pic(C)$ is equal to the quotient group of the group of divisors on $C$ modulo the group of principal divisors on $C$. We always write $\mathcal{O}_C(D)$ for its equivalence class in $\Pic(C)$ unless there might be ambiguity.

The degree map
\begin{align*}
    \deg:\Pic(C)&\rightarrow\mathbb{Z}\\
    \mathcal{O}_C(D)&\mapsto\deg D
\end{align*}
is a surjective group homomorphism. We denote by $\Pic^0(C)$ the kernel of $\deg$. Consider the map
\begin{align*}
    C(k)&\rightarrow\Pic^0(C)\\
    x&\mapsto\mathcal{O}_C(x-O).
\end{align*}
By \cite[Chapter III, Proposition 3.4(d); Chapter X, Theorem 3.8(c)]{Sil2009}, it is bijective, and hence induces a group structure on $C(k)$. The above arguments give a decomposition of $\Pic(C)$:
\begin{equation}\label{Picard group}
    \Pic(C)\cong\mathbb{Z}\oplus C(k).
\end{equation}
More precisely, if $\mathcal{F}$ is a line bundle over $C$ with degree $d$, then we must have
\[\mathcal{F}\cong\mathcal{O}_C(x-O)\otimes\mathcal{O}_C(dO)\]
for some $x\in C(k)$. By the finiteness of $k$, $C(k)$ is a finite group, and hence $\Pic(C)$ is finitely generated.

To provide the structure of the Grothendieck group $K_0(\Coh C)$, we have to introduce the notion of the determinant of a coherent sheaf. Let $\mathcal{F}$ be a vector bundle over $C$ of rank $r$. We define its \emph{determinant} to be the exterior power
\[\det\mathcal{F}=\wedge^r\mathcal{F}.\]
It is a line bundle over $C$. For a torsion sheaf $\mathcal{T}$ over $C$, we set the divisor
\[D_\mathcal{T}=\sum_{x\in\supp\mathcal{T}}\ell(\mathcal{T},x)x.\]
For a general coherent sheaf $\mathcal{F}$ over $C$, we set
\[\det\mathcal{F}=\det(\mathcal{F}/\mathcal{F}_\tor)\otimes\mathcal{O}_C(D_{\mathcal{F}_\tor}).\]
It is also a line bundle over $C$. We can extend the degree function to arbitrary coherent sheaves:
\[\deg\mathcal{F}:=\deg(\det\mathcal{F})=\deg(\det(\mathcal{F}/\mathcal{F}_\tor))+\sum_{x\in\supp\mathcal{F}_\tor}\ell(\mathcal{F}_\tor,x)[k_x:k].\]
The following proposition provides a characterization of $K_0(\Coh C)$.
\begin{proposition}
    The map
    \begin{equation}\label{K0(COh C)}
        \begin{split}
        K_0(\Coh C)&\rightarrow\mathbb{Z}\oplus\Pic(C)\\
        \udv\mathcal{F}&\mapsto(\rk\mathcal{F},\det\mathcal{F})
    \end{split}
    \end{equation}
    is a well-defined group isomorphism. In particular, 
    \[K_0(\Coh C)\cong\mathbb{Z}\oplus\Pic(C)\cong\mathbb{Z}^2\oplus C(k)\]
    is finitely generated, and the \emph{numerical class function}
    \begin{equation}\label{numerical class function}
        \begin{split}
        \cl:K_0(\Coh C)&\rightarrow\mathbb{Z}^2\\
        \udv\mathcal{F}&\mapsto(\rk\mathcal{F},\deg\mathcal{F})
    \end{split}
    \end{equation}
    is a well-defined surjective group homomorphism.
\end{proposition}
\begin{proof}
    By \cite[Chapter II, Theorem 8.2 and Proposition 8.2.1]{Wei2013}, the map (\ref{K0(COh C)}) is a well-defined group isomorphism, and the rest of the proposition follows from (\ref{Picard group}) and the surjectivity of $\deg$.
\end{proof}
\subsection{Slopes and the HN filtrations}
Let $0\neq\mathcal{F}\in\Coh C$. We define the \emph{slope} of $\mathcal{F}$ to be
\[\mu(\mathcal{F})=\begin{cases}
    \dfrac{\deg\mathcal{F}}{\rk\mathcal{F}}&,\rk\mathcal{F}>0,\\
    \infty&,\rk\mathcal{F}=0.
\end{cases}\]
\begin{definition}
    Let $0\neq\mathcal{F}\in\Coh C$. We say that $\mathcal{F}$ is \emph{semistable} if $\mu(\mathcal{E})\leqslant\mu(\mathcal{F})$ for any nonzero proper subsheaf $\mathcal{E}$ of $\mathcal{F}$. We say that $\mathcal{F}$ is \emph{stable} if $\mu(\mathcal{E})<\mu(\mathcal{F})$ for any nonzero proper subsheaf $\mathcal{E}$ of $\mathcal{F}$.
\end{definition}
\begin{example}
    \begin{enumerate}
        \item Line bundles are stable. In fact, a nonzero proper subsheaf of a line bundle $\mathcal{O}_C(D)$ must have the form $\mathcal{O}_C(D-D')$, where $D'$ is a nonzero effective divisor, and
    \[\mu(\mathcal{O}_C(D-D'))=\deg(D-D')<\deg D=\mu(\mathcal{O}_C(D)).\]
    \item A direct sum of more than one line bundle of the same degree is semistable but cannot be stable, since each line bundle is a nonzero proper subsheaf with the same slope.
    \item Torsion sheaves are semistable, since a nonzero subsheaf of a torsion sheaf must be torsion, and hence has slope $+\infty$. A torsion sheaf is stable if and only if it is simple; that is, it is a skyscraper sheaf.
    \end{enumerate}
\end{example}
\begin{theorem}[{\cite[Proposition 1.3.9]{HN1975}\cite[\S 1.1]{BS2012}}]
    Let $0\neq\mathcal{F}\in\Coh C$. Then $\mathcal{F}$ admits a unique \emph{Harder--Narasimhan filtration}
    \[0=\mathcal{F}_0\subsetneqq\mathcal{F}_1\subsetneqq\dots\subsetneqq\mathcal{F}_m=\mathcal{F},\]
    such that each quotient $\mathcal{G}_i:=\mathcal{F}_i/\mathcal{F}_{i-1},i=1,\dots,m$ is semistable, with
    \[\mu(\mathcal{G}_1)>\mu(\mathcal{G}_2)>\dots>\mu(\mathcal{G}_m).\]
    In particular, $\mathcal{G}_1$ is torsion if and only if $\mathcal{F}$ is not a vector bundle. We set
    \[\mu_-(\mathcal{F}):=\mu(\mathcal{G}_m),\mu_+(\mathcal{F}):=\mu(\mathcal{G}_1).\]
\end{theorem}
We have the following lemma of orthogonality.
\begin{lemma}\label{orthogonality}
    Let $0\neq\mathcal{F},\mathcal{G}\in\Coh C$. If $\mu_-(\mathcal{F})>\mu_+(\mathcal{G})$, then
    \[\Hom_{\Coh C}(\mathcal{F},\mathcal{G})=0.\]
\end{lemma}
For $\mu\in\mathbb{R}\cup\{\infty\}$, we set the subcategory of $\Coh C$:
\[\mathcal{C}_\mu=\{0\neq\mathcal{F}\in\Coh C|\mathcal{F}\text{ is semistable and }\mu(\mathcal{F})=\mu\}\cup\{0\}.\]
Note that when $\mu\in\mathbb{R}-\mathbb{Q}$, the category $\mathcal{C}_\mu$ must be $\{0\}$.
For any $\mu\in\mathbb{Q}\cup\{\infty\}$, $\mathcal{C}_\mu$ is an abelian category closed under extension. We call $\mathcal{C}_\mu$ the \emph{slope category} of slope $\mu$. The simple objects in $\mathcal{C}_\mu$ are just the stable objects of slope $\mu$. Note that $\mathcal{C}_{\infty}$ coincides with the subcategory of torsion sheaves.

The following property is very important.
\begin{proposition}[{\cite[\S 2]{BS2012}}]\label{semistable sheaf finite}
    For any $a\in\mathbb{N},b\in\mathbb{Z}$, the set
    \[\{\mathcal{F}\in\Iso(\Coh C)\text{ semistable}|\rk\mathcal{F}=a,\deg\mathcal{F}=b\}\]
    is finite.
\end{proposition}
\subsection{The stable tubes}
Let $\mathcal{S}\in\Coh C$ be a stable object. We define the \emph{stable tube} $T(\mathcal{S})$ to be the subcategory of $\mathcal{C}_{\mu(\mathcal{S})}$ whose objects have composition factors all isomorphic to $\mathcal{S}$.

Let us first focus on the category $\mathcal{C}_\infty$ of torsion sheaves. Recall that a torsion sheaf $\mathcal{F}$ has a decomposition
\[\mathcal{F}=\bigoplus_{x\in\supp\mathcal{F}}\bigoplus_{j=1}^{r_x}\mathcal{O}_C/\mathcal{I}_x^{n_{x,j}},n_{x,j}\in\mathbb{Z}_+.\]
Thus an indecomposable torsion sheaf over $C$ must be supported on a single point $x$. Note that the category $\mathcal{C}_{\infty,x}$ of torsion sheaves supported on $x$ is equivalent to the category of finite-length $\mathcal{O}_{C,x}$-modules, and hence equivalent to the category of finite-dimensional nilpotent $k_x[t]$-modules, where $t$ is a formal indeterminate. The equivalence is given by
\[\mathcal{O}_C/\mathcal{I}_x^a\mapsto k_x[t]/(t^a).\]
We will always identify the two categories via this equivalence, and we will write
\[J_a=k_x[t]/(t^a)\in\mathcal{C}_{\infty,x},a\in\mathbb{N},\]
where $J_0=0$. For $a>0$, the following composition series gives all subsheaves of $J_a$:
\[0=t^aJ_a\subsetneqq t^{a-1}J_a\subsetneqq\dots\subsetneqq tJ_a\subsetneqq J_a,\]
with each composition factor being $J_1=\underline{k}_x$. So we have $J_a\in T(\underline{k}_x)$ for any $a\in\mathbb{Z}_+$ and hence $\mathcal{C}_{\infty,x}=T(\underline{k}_x)$. Thus we have a decomposition of $\mathcal{C}_\infty$:
\begin{equation}
    \mathcal{C}_\infty=\bigoplus_{x\in |C|}\mathcal{C}_{\infty,x}=\bigoplus_{x\in |C|}T(\underline{k}_x).
\end{equation}
Moreover, there is no nonzero morphism or nontrivial extension between an object in $T(\underline{k}_x)$ and an object in $T(\underline{k}_y)$ for any $x\neq y\in |C|$.

For other slope categories, we have
\begin{lemma}[{\cite[Theorem 1.1(iii)]{BS2012}}]
    For any $\mu\in\mathbb{Q}$, there exists an equivalence of categories:
    \[\varepsilon_\mu:\mathcal{C}_\mu\rightarrow\mathcal{C}_\infty.\]
\end{lemma}
For any stable object $\mathcal{S}\in\mathcal{C}_\mu$, we set
\[K_{\mathcal{S}}=\End_{\Coh C}\mathcal{S}.\]
It is a finite-dimensional algebra over $k$, and hence finite as a set. Since $\mathcal{S}$ is simple in $\mathcal{C}_\mu$, $K_{\mathcal{S}}$ is a division ring by Schur's lemma. Thus, $K_{\mathcal{S}}$ is a field and $[K_{\mathcal{S}}:k]<+\infty$. The stable tube $T(\mathcal{S})$ is equivalent to the category of finite-dimensional nilpotent $K_\mathcal{S}[t]$-modules. The object in $T(\mathcal{S})$ corresponding to $K_\mathcal{S}[t]/(t^a)$ via the equivalence is denoted by $J_a(\mathcal{S}),a\in\mathbb{N}$. When $a>0$, its composition series exhibits all of its subobjects:
\begin{equation}\label{composition series}
    0=t^aJ_a(\mathcal{S})\subsetneqq t^{a-1}J_a(\mathcal{S})\subsetneqq\dots\subsetneqq tJ_a(\mathcal{S})\subsetneqq J_a(\mathcal{S}).
\end{equation}
\begin{lemma}\label{Tabl}
    Let $a,b\in\mathbb{Z}_+$. Up to automorphisms of both objects, a nonzero morphism $J_a(\mathcal{S})\rightarrow J_b(\mathcal{S})$ can be written as
    \begin{equation}\label{fl}
        \begin{split}
            f_l:K_\mathcal{S}[t]/(t^a)&\rightarrow K_\mathcal{S}[t]/(t^b),\\
            1&\mapsto t^{b-l},
        \end{split}
    \end{equation}
    where $1\leqslant l\leqslant\min(a,b)$. Moreover,
    \[\Ker f_l\cong J_{a-l}(\mathcal{S}),\im f_l\cong J_l(\mathcal{S}),\Coker f_l\cong J_{b-l}(\mathcal{S}).\]
\end{lemma}
\begin{proof}
    Let $f:K_\mathcal{S}[t]/(t^a)\rightarrow K_\mathcal{S}[t]/(t^b)$ be a nonzero $K_\mathcal{S}[t]$-module homomorphism. We can write
    \[f(1)=t^{b-l}h(t),1\leqslant l\leqslant b,\]
    where $h(0)\neq0$ and hence $h(t)$ is a unit in $K_\mathcal{S}[t]/(t^b)$. Since $f(0)=f(t^a)$ has to be zero, we obtain $t^{a+b-l}=0$ in $K_\mathcal{S}[t]/(t^b)$, so $l\leqslant a$. Thus $f$ is the composition of $f_l$ and the automorphism on $K_\mathcal{S}[t]/(t^b)$ defined by multiplication with $h(t)$. One can directly compute the kernel, image and cokernel of $f_l$.
\end{proof}
\subsection{Morphism categories and elliptic triples}
Let $\mathcal{C}$ be a category. The \emph{morphism category} $\Mor(\mathcal{C})$ is defined to be the category consisting of the following data:
\begin{itemize}
    \item An object in $\Mor(\mathcal{C})$ is a triple $(A,B,f)$, where $A,B\in\mathcal{C}$ and $f\in\Hom_\mathcal{C}(A,B)$. We always write such a triple as $\left(A\xrightarrow{f}B\right)$.
    \item A morphism from $\left(A\xrightarrow{f}B\right)$ to $\left(A'\xrightarrow{f'}B'\right)$ is a pair $(g_1,g_2)$, such that $g_1\in\Hom_\mathcal{C}(A,A'),g_2\in\Hom_\mathcal{C}(B,B')$ and $g_2\circ f=f'\circ g_1$, or equivalently, the diagram
    \begin{center}
        \begin{tikzcd}
A \arrow[r, "f"] \arrow[d, "g_1"'] & B \arrow[d, "g_2"] \\
A' \arrow[r, "f'"]                 & B'                
\end{tikzcd}
    \end{center}
    commutes.
\end{itemize}

The following property of morphism categories is a special case of \cite[Theorem 1]{JM2010}.
\begin{proposition}\label{morphism category}
    Let $\mathcal{A}$ be an abelian category. Then the morphism category $\Mor(\mathcal{A})$ is also abelian, with the kernel and cokernel of a morphism $(g_1,g_2):\left(A\xrightarrow{f}B\right)\rightarrow\left(A'\xrightarrow{f'}B'\right)$ defined to be
    \begin{align*}
        \Ker(g_1,g_2)&=\left(\Ker g_1\xrightarrow{f|_{\Ker g_1}}\Ker g_2\right),\\
        \Coker(g_1,g_2)&=\left(\Coker g_1\xrightarrow{\overline{f'}}\Coker g_2\right)
    \end{align*}
    together with the canonical morphisms.
\end{proposition}
Using Proposition \ref{morphism category}, one can easily show that the exactness in $\Mor(\mathcal{A})$ can be described on the objects on each vertex. We also refer the reader to \cite[Proposition 2.4.1]{ZW2018}.

Now take $\mathcal{A}=\Coh C$. We set $\TCoh C=\Mor(\Coh C)$. By Proposition \ref{morphism category}, $\TCoh C$ is an abelian category.

Given an object $\mathcal{F}\in\Coh C$, we set the following three triples:
\[A_\mathcal{F}=(\mathcal{F}\rightarrow0),B_\mathcal{F}=(0\rightarrow\mathcal{F}),D_\mathcal{F}=(\mathcal{F}\xrightarrow{\id_\mathcal{F}}\mathcal{F}).\]
For any $T=(\mathcal{F}\xrightarrow{f}\mathcal{G})\in\TCoh C$, we have a short exact sequence
\begin{equation}\label{decomposition}
    0\rightarrow B_\mathcal{G}\rightarrow T\rightarrow A_\mathcal{F}\rightarrow0.
\end{equation}
Thus the Grothendieck group of $\TCoh C$ is
\[K_0(\TCoh C)=K_0(\Coh C)\oplus K_0(\Coh C)\cong\mathbb{Z}^4\oplus C(k)^2.\]
When $|C(k)|>1$, one can easily see that $K_0(\TCoh C)$ is not free, and hence $\TCoh C$ cannot be of finite length. We fix a \emph{numerical class function}:
\begin{equation}\label{cl}
    \begin{split}
        \cl: K_0(\TCoh C)&\rightarrow\mathbb{Z}^4\\
        \udv(\mathcal{F}\xrightarrow{f}\mathcal{G})&\mapsto(\rk\mathcal{F},\deg\mathcal{F},\rk\mathcal{G},\deg\mathcal{G}).
    \end{split}
\end{equation}
\begin{proposition}[{\cite[Proposition 2.7.24]{Rom2018}}]\label{homdim2}
    The homological dimension of $\TCoh C$ is $2$.
\end{proposition}

\section{Bridgeland stability conditions on elliptic triples}
\subsection{Bridgeland stability conditions}
We give a short review of the theory of Bridgeland stability conditions. We refer the reader to \cite{Bri2007} and \cite[\S 2.2-2.3]{GX2026} for details.
\begin{definition}
    Let $\mathcal{D}$ be a triangulated category. A \emph{Bridgeland stability condition} on $\mathcal{D}$ is a pair $\sigma=(Z,\mathcal{P})$, where
    \begin{itemize}
        \item $Z:K_0(\mathcal{D})\rightarrow\mathbb{C}$ is a group homomorphism, called the \emph{central charge};
        \item $\mathcal{P}=\{\mathcal{P}(\phi)|\phi\in\mathbb{R}\}$ is a family of full additive extension-closed subcategories of $\mathcal{D}$,
    \end{itemize}
    satisfying the following conditions:
    \begin{enumerate}
        \item If $0\neq M\in\mathcal{P}(\phi)$, then $Z(M)=m_Me^{i\pi\phi}$ for some $m_M\in\mathbb{R}_+$.
        \item $\mathcal{P}(\phi+1)=\mathcal{P}(\phi)[1]$ for any $\phi\in\mathbb{R}$.
        \item For $A_i\in\mathcal{P}(\phi_i),i=1,2$ with $\phi_1>\phi_2$, we have $\Hom_\mathcal{D}(A_1,A_2)=0$.
        \item (HN filtration) For any $0\neq E\in\mathcal{D}$, there is a sequence of triangles
        \begin{center}
            \begin{tikzcd}[sep=tiny]
0=E_0 \arrow[rr] &                        & E_1 \arrow[rr] \arrow[ld] &                        & E_2 \arrow[r] \arrow[ld] & \cdots \arrow[r] & E_{n-1} \arrow[rr] &                        & E_n=E \arrow[ld] \\
                 & A_1 \arrow[lu, dashed] &                           & A_2 \arrow[lu, dashed] &                          &                  &                    & A_n \arrow[lu, dashed] &                 
\end{tikzcd}
        \end{center}
        such that each $A_i\neq0$ belongs to $\mathcal{P}(\phi_i)$ for some $\phi_i$, $1\leqslant i\leqslant n$, with $\phi_1>\cdots>\phi_n$. This sequence is called a \emph{Harder-Narasimhan filtration}, or an \emph{HN filtration} for short. It is unique up to isomorphism by \cite{Bri2007}. Thus we can define
        \[\phi^+(E)=\phi_1,\phi^-(E)=\phi_n,\Delta(E)=\phi_1-\phi_n\]
        and
        \[m(E):=\sum_{j=1}^n|Z(A_j)|.\]
    \end{enumerate}
\end{definition}

For an interval $I\subseteq\mathbb{R}$, let $\mathcal{P}(I)$ consist of zero and the objects whose HN factors have phases in $I$, that is,
\[\mathcal{P}(I)=\{0\neq E\in\mathcal{D}|\phi^{\pm}(E)\in I\}\cup\{0\}.\]
It is extension-closed. The subcategory $\mathcal{P}(0,1]$ is called the \emph{heart} of $\sigma$. It is an abelian category by \cite[\S 3]{Bri2007}. For any $\phi\in\mathbb{R}$, the subcategory $\mathcal{P}(\phi,\phi+1]$ is also abelian and called a \emph{rotated heart}. Bridgeland proved the following important result:
\begin{lemma}[{\cite[Proposition 5.3]{Bri2007}}]\label{heart}
    To give a stability condition on a triangulated category $\mathcal{D}$ is equivalent to giving a bounded t-structure on $\mathcal{D}$ and a stability function having the HN property on its heart; see \cite[Definitions 2.1-2.3]{Bri2007}.
\end{lemma}

By \cite[Lemma 5.2]{Bri2007}, the subcategories $\mathcal{P}(\phi),\phi\in\mathbb{R}$ are abelian. We call $\mathcal{P}(\phi)$ the \emph{phase category} with phase $\phi$. It will play a major role in our study. We say an object $0\neq E\in\mathcal{D}$ is \emph{$\sigma$-semistable} if $E$ belongs to some phase category $\mathcal{P}(\phi)$ for some phase $\phi\in\mathbb{R}$, and we say it is \emph{$\sigma$-stable} if it is simple in the abelian phase category $\mathcal{P}(\phi)$. For a $\sigma$-semistable object $E\in\mathcal{P}(\phi)$, we often write $\phi(E)=\phi$.
\begin{proposition}[{\cite[Proposition 2.3]{GX2026}}]\label{phase truncation}
    Let $E\in\mathcal{D}$ and $c\in\mathbb{R}$. There is a unique (up to unique isomorphism of triangles inducing the identity on $E$) \emph{phase truncation triangle}
    \[E_{>c}\rightarrow E\rightarrow E_{\leqslant c}\rightarrow E_{>c}[1],\]
    where $E_{>c}\in\mathcal{P}(c,+\infty)$ and $E_{\leqslant c}\in\mathcal{P}(-\infty,c]$.
\end{proposition}
\begin{definition}[{\cite[Definition 5.7]{Bri2007}}]
    We say a stability condition $\sigma=(Z,\mathcal{P})$ on $\mathcal{D}$ is \emph{locally finite} if there exists $0<\eta<\dfrac12$ such that the quasi-abelian category $\mathcal{P}(\phi-\eta,\phi+\eta)$ is of finite length for any $\phi\in\mathbb{R}$.
\end{definition}
Let $\Stab(\mathcal{D})$ be the set of locally finite stability conditions on $\mathcal{D}$. We introduce the generalized metric, that is, a metric that allows infinite distance, on $\Stab(\mathcal{D})$ defined by Bridgeland:
\begin{lemma}[{\cite[Proposition 8.1]{Bri2007}}]
    Let $\sigma,\tau\in\Stab(\mathcal{D})$. The function
    \[d(\sigma,\tau)=\sup_{0\neq E\in\mathcal{D}}\max\left\{|\phi_\sigma^-(E)-\phi_\tau^-(E)|,|\phi_\sigma^+(E)-\phi_\tau^+(E)|,\left|\log\dfrac{m_\sigma(E)}{m_\tau(E)}\right|\right\}\]
    defines a generalized metric on $\Stab(\mathcal{D})$.
\end{lemma}
The generalized metric $d$ induces a topology on $\Stab(\mathcal{D})$.
Let us fix a group epimorphism
\[\cl:K_0(\mathcal{D})\rightarrow\Gamma,\]
where $\Gamma$ is a finite rank lattice. For an object $E\in\mathcal{D}$, we always write $\cl E=\cl(\udv E)$ for short.
\begin{definition}
    A stability condition $\sigma=(Z,\mathcal{P})$ is called \emph{numerical} with respect to $\cl$ if $Z$ factors through $\cl$. The subspace of $\Stab(\mathcal{D})$ consisting of numerical stability conditions is denoted by $\Stab_\cl(\mathcal{D})$.
\end{definition}
\begin{proposition}[{\cite[Corollary 1.3]{Bri2007}}]
    The subspace $\Stab_\cl(\mathcal{D})$ is a complex manifold.
\end{proposition}
Fix a norm $\|\cdot\|$ on the complex vector space $\Gamma_{\mathbb{C}}=\Gamma\otimes_{\mathbb{Z}}\mathbb{C}$. We introduce the following concepts:
\begin{definition} Let $\sigma=(Z,\mathcal{P})\in\Stab_\cl(\mathcal{D})$.
    \begin{enumerate}
        \item $\sigma$ is said to have the \emph{support property} if there exists a $c>0$ such that
        \[|Z(\cl E)|\geqslant c\|\cl E\|\]
        for any $\sigma$-semistable object $E$.
        \item $\sigma$ is said to have the \emph{finite counting condition} if for any $\phi\in\mathbb{R}$ and $\gamma\in\Gamma$, the set
    \[\{E\in\Iso(\mathcal{D})|E\in\mathcal{P}(\phi),\cl E=\gamma\}\]
    is finite.
    \end{enumerate}
\end{definition}
When $\mathcal{D}$ is $k$-linear and Hom-finite, we have the following properties.
\begin{proposition}[{\cite[Proposition 2.10]{GX2026}}]\label{mass and phase bounds finite}
    Let $\sigma=(Z,\mathcal{P})\in\Stab_\cl(\mathcal{D})$ satisfy the support property and the finite counting condition. For any $R>0$ and any bounded segment $I\subseteq\mathbb{R}$, the set
    \[\{0\neq E\in\Iso(\mathcal{P}(I))|m(E)\leqslant R\}\]
    is finite.
\end{proposition}
\begin{corollary}[{\cite[Corollary 2.11]{GX2026}}]\label{finite distance}
    Let $\sigma=(Z_\sigma,\mathcal{P}_\sigma)\in\Stab_\cl(\mathcal{D})$ satisfy the support property and the finite counting condition. If $\tau=(Z_\tau,\mathcal{P}_\tau)\in\Stab_\cl(\mathcal{D})$ with $d(\sigma,\tau)<+\infty$, then $\tau$ satisfies the support property and the finite counting condition.
\end{corollary}
\subsection{Stability conditions for elliptic triples}
Let $\mathcal{T}_C=D^b(\TCoh C)$ be the bounded derived category of the elliptic triple category $\TCoh C$ and fix $\cl$ as (\ref{cl}) equipped with the $1$-norm $\|\cdot\|_1$. We now construct a type of Bridgeland stability conditions on $\mathcal{T}_C$. Our construction follows \cite[Definition 2.1]{BGG2004}.
\begin{lemma}[{\cite[Lemma 5.3.9; Proposition 5.3.10]{Dal2007}}]\label{ABC}
    Let $Z:\mathbb{Z}^4=\cl(K_0(\TCoh C))\rightarrow\mathbb{C}$ be a group homomorphism. Then $Z$ is a stability function on $\TCoh C$ with HN property inducing a locally finite stability condition on $\mathcal{T}_C$, in the sense of Lemma \ref{heart}, if and only if $Z$ has the form
    \begin{equation}\label{general stability function}
        Z(r_1,d_1,r_2,d_2)=-(A_1d_1+A_2d_2)+(B_1r_1+B_2r_2)+(C_1r_1+C_2r_2)i,
    \end{equation}
    where $A_1,A_2,B_1,B_2,C_1,C_2\in\mathbb{R}$ with $A_1,A_2,C_1,C_2>0$.
\end{lemma}

Let $\alpha\in\mathbb{R}$. Consider the function
\begin{equation}
    \begin{split}
        Z_\alpha:\mathbb{Z}^4=\cl(K_0(\TCoh C))&\rightarrow\mathbb{C}\\
        (r_1,d_1,r_2,d_2)&\mapsto-(d_1+d_2+\alpha r_1)+(r_1+r_2)i.
    \end{split}
\end{equation}
Note that $Z_\alpha$ is obtained by taking $A_1=A_2=C_1=C_2=1,B_1=-\alpha,B_2=0$ in (\ref{general stability function}). So by Lemma \ref{ABC}, $Z_\alpha$ is a stability function on $\TCoh C$ with HN property, and the induced stability condition $\sigma_\alpha=(Z_\alpha,\mathcal{P}_\alpha)$ on $\mathcal{T}_C$ is locally finite. For $0\neq T=(\mathcal{F}\xrightarrow{f}\mathcal{G})\in\TCoh C$, let $\cl T=(r_1,d_1,r_2,d_2)$. We define the \emph{$Z_\alpha$-slope} $\nu_\alpha(T)\in\mathbb{R}\cup\{\infty\}$ of $T$ to be
\[\nu_\alpha(T)=\begin{cases}
    \dfrac{d_1+d_2+\alpha r_1}{r_1+r_2}=-\cot(\arg Z_\alpha(T)),&r_1+r_2>0,\\
    \infty,&r_1+r_2=0,
\end{cases}\]
where $\arg$ is the branch of the argument taking values in $(0,\pi]$.
$T$ is torsion, that is, $\mathcal{F}$ and $\mathcal{G}$ are both torsion, if and only if $\nu_\alpha(T)=+\infty$. One can easily see that a torsion triple must be $Z_\alpha$-semistable.
\begin{lemma}\label{slope bounds}
    Let $T=(\mathcal{F}\xrightarrow{f}\mathcal{G})\in\TCoh C$ be $Z_\alpha$-semistable, with $r_1+r_2>0$. Set $\nu=\nu_\alpha(T)$.
    \begin{enumerate}
        \item $\mathcal{F}$ and $\mathcal{G}$ are both torsion-free.
        \item If $\mathcal{F}\neq0$, then $\mu_-(\mathcal{F})\geqslant\nu-\alpha$. If $\mathcal{G}\neq0$, then $\mu_+(\mathcal{G})\leqslant\nu$.
        \item If $\mathcal{F}\neq0$, then $\mu_+(\mathcal{F})\leqslant\max(\nu,\nu-\alpha)$. If $\mathcal{G}\neq0$, then $\mu_-(\mathcal{G})\geqslant\min(\nu,\nu-\alpha)$.
        \item If $\alpha\leqslant0$, then $\mathcal{F}\in\mathcal{C}_{\nu-\alpha},\mathcal{G}\in\mathcal{C}_\nu$. If $\alpha<0$, then $f=0$.
    \end{enumerate}
\end{lemma}
\begin{proof}
    (1) The triple $T_\tor=(\mathcal{F}_\tor\xrightarrow{f|_{\mathcal{F}_\tor}}\mathcal{G}_\tor)$ is a subobject of $T$. If $T_\tor\neq0$, then $\nu_\alpha(T_\tor)=+\infty>\nu_\alpha(T)$, which contradicts the assumption that $T$ is $Z_\alpha$-semistable. So $T_\tor=0$ and hence $\mathcal{F}$ and $\mathcal{G}$ are both torsion-free.

    (2) When $\mathcal{F}\neq0$, let $\mathcal{F}''$ be a nonzero quotient sheaf of $\mathcal{F}$. Then $A_{\mathcal{F}''}$ is a quotient object of $T$ in $\TCoh C$. By the $Z_\alpha$-semistability of $T$, we have
    \[\nu=\nu_\alpha(T)\leqslant\nu_\alpha(A_{\mathcal{F}''})=\mu(\mathcal{F}'')+\alpha.\]
    Now take $\mathcal{F}''$ to be the last HN factor of $\mathcal{F}$. Then $\mu(\mathcal{F}'')=\mu_-(\mathcal{F})$, and hence $\mu_-(\mathcal{F})\geqslant\nu-\alpha$.

    When $\mathcal{G}\neq0$, let $\mathcal{G}'$ be a nonzero subsheaf of $\mathcal{G}$. Then $B_{\mathcal{G}'}$ is a subobject of $T$ in $\TCoh C$. By the $Z_\alpha$-semistability of $T$, we have
    \[\nu=\nu_\alpha(T)\geqslant\nu_\alpha(B_{\mathcal{G}'})=\mu(\mathcal{G}').\]
    Now take $\mathcal{G}'$ to be the first HN factor. Then $\mu(\mathcal{G}')=\mu_+(\mathcal{G})$, and hence $\mu_+(\mathcal{G})\leqslant\nu$.

    (3) When $\mathcal{F}\neq0$, let $\mathcal{F}'$ be the first HN factor of $\mathcal{F}$. If $\mu_+(\mathcal{F})=\mu(\mathcal{F}')>\nu$, then $\Hom_{\Coh C}(\mathcal{F'},\mathcal{G})=0$. This is clear if $\mathcal{G}=0$; otherwise, it follows from (2) and Lemma \ref{orthogonality}. In particular, $f|_{\mathcal{F}'}=0$. So $A_{\mathcal{F}'}$ is a subobject of $T$ in $\TCoh C$. Therefore,
    \[\nu=\nu_\alpha(T)\geqslant\nu_\alpha(A_{\mathcal{F}'})=\mu(\mathcal{F}')+\alpha=\mu_+(\mathcal{F})+\alpha,\]
    that is, $\mu_+(\mathcal{F})\leqslant\nu-\alpha$. Hence, we have $\mu_+(\mathcal{F})\leqslant\max(\nu,\nu-\alpha)$.

    When $\mathcal{G}\neq0$, let $\mathcal{G}''$ be the last HN factor of $\mathcal{G}$. If $\mu_-(\mathcal{G})=\mu(\mathcal{G}'')<\nu-\alpha$, then $\Hom_{\Coh C}(\mathcal{F},\mathcal{G}'')=0$. This is clear if $\mathcal{F}=0$; otherwise, it follows from (2) and Lemma \ref{orthogonality}. In particular, the composite morphism $\mathcal{F}\xrightarrow{f}\mathcal{G}\twoheadrightarrow\mathcal{G''}$ is zero. So $B_{\mathcal{G}''}$ is a quotient object of $T$ in $\TCoh C$. Therefore,
    \[\nu=\nu_\alpha(T)\leqslant\nu_\alpha(B_{\mathcal{G}''})=\mu(\mathcal{G}'')=\mu_-(\mathcal{G}),\]
    that is, $\mu_-(\mathcal{G})\geqslant\nu$. Hence, we have $\mu_-(\mathcal{G})\geqslant\min(\nu,\nu-\alpha)$.

    (4) When $\alpha\leqslant0$, we have $\nu-\alpha\geqslant\nu$. Then the results follow by (2), (3) and Lemma \ref{orthogonality}.
\end{proof}
For $0\neq T\in\TCoh C$, we set
\[\theta_\alpha(T)=\dfrac{1}{\pi}\arg Z_\alpha(\cl T)\in(0,1].\]
\begin{lemma}\label{phase control}
    For any $0\neq T\in\TCoh C$ and $\alpha\neq\beta\in\mathbb{R}$, we have
    \[|\theta_\alpha(T)-\theta_\beta(T)|\leqslant\dfrac{|\alpha-\beta|}{\pi}.\]
\end{lemma}
\begin{proof}
    Let $T=(\mathcal{F}\xrightarrow{f}\mathcal{G})$ and $\cl T=(r_1,d_1,r_2,d_2)$. When $r_1=r_2=0$, that is, when $T$ is torsion, $\theta_\alpha(T)=\theta_\beta(T)=1$, so the inequality holds trivially. When $r_1+r_2>0$, the function
    \[\arg Z_t(T)=\dfrac{\pi}{2}+\arctan\dfrac{d_1+d_2+tr_1}{r_1+r_2}\]
    is smooth in $t\in\mathbb{R}$, with
    \[\left|\dfrac{\diff}{\diff t}\arg Z_t(T)\right|=\dfrac{r_1(r_1+r_2)}{(d_1+d_2+tr_1)^2+(r_1+r_2)^2}\leqslant1,t\in\mathbb{R}.\]
    So by Lagrange's mean value theorem, we have
    \[\left|\dfrac{\arg Z_\alpha(T)-\arg Z_\beta(T)}{\alpha-\beta}\right|\leqslant1,\]
    and the inequality we want to prove follows immediately.
\end{proof}

\begin{theorem}\label{support}
    Let $\alpha\in\mathbb{R}$. The stability condition $\sigma_\alpha=(Z_\alpha,\mathcal{P}_\alpha)$ on $\mathcal{T}_C$ satisfies the support property, with a support constant
    \[c_\alpha=\dfrac{1}{\sqrt{1+(1+|\alpha|)^2}}.\]
\end{theorem}
\begin{proof}
    We are going to show that for any $T\in\mathcal{T}_C$ that is $\sigma_\alpha$-semistable, we have
    \[|Z_\alpha(\cl T)|\geqslant c_\alpha\|\cl T\|_1.\]
    Without loss of generality, we may assume that $T\in\mathcal{P}_\alpha(0,1]=\TCoh C$, after shifting $T$ if necessary. Let $\cl T=(r_1,d_1,r_2,d_2),\nu=\nu_\alpha(T)$. When $r_1+r_2=0$, $T$ is torsion and hence $r_1=r_2=0,d_1,d_2\geqslant0$ and $d_1+d_2>0$. So $Z_\alpha(\cl T)=-(d_1+d_2)$. Since $c_\alpha<1$, we have
    \[|Z_\alpha(\cl T)|=d_1+d_2=\|\cl T\|_1\geqslant c_\alpha\|\cl T\|_1.\]
    When $r_1+r_2>0$, we have $\nu\in\mathbb{R}$. Write $T=(\mathcal{F}\xrightarrow{f}\mathcal{G})$. When $\mathcal{F}\neq0$, let $\mathcal{E}_1,\dots,\mathcal{E}_n$ be the HN factors of $\mathcal{F}$. By Lemma \ref{slope bounds}(2) and (3), we have
    \[\nu-\alpha\leqslant\mu(\mathcal{E}_i)\leqslant\max(\nu,\nu-\alpha).\]
    Hence $|\mu(\mathcal{E}_i)|\leqslant|\nu|+|\alpha|,i=1,\dots,n$. Thus we have
    \begin{align*}
        |d_1|&=\left|\sum_{i=1}^n\deg\mathcal{E}_i\right|\\
        &=\left|\sum_{i=1}^n\mu(\mathcal{E}_i)\rk\mathcal{E}_i\right|\\
        &\leqslant\sum_{i=1}^n|\mu(\mathcal{E}_i)|\rk\mathcal{E}_i\\
        &\leqslant\sum_{i=1}^n(|\nu|+|\alpha|)\rk\mathcal{E}_i\\
        &=(|\nu|+|\alpha|)\sum_{i=1}^n\rk\mathcal{E}_i\\
        &=(|\nu|+|\alpha|)r_1.
    \end{align*}
When $\mathcal{F}=0$, we have $d_1=0=r_1$. Thus we have an inequality
\[|d_1|\leqslant(|\nu|+|\alpha|)r_1.\]
Similarly, one can show that
\[|d_2|\leqslant(|\nu|+|\alpha|)r_2.\]
Thus,
\begin{align*}
    \|\cl T\|_1&=r_1+r_2+|d_1|+|d_2|\\
    &\leqslant(1+|\nu|+|\alpha|)(r_1+r_2)\\
    &=(1+|\alpha|)(r_1+r_2)+|d_1+d_2+\alpha r_1|\\
    (\text{Cauchy-Schwarz})&\leqslant\sqrt{(1+|\alpha|)^2+1}\sqrt{(r_1+r_2)^2+(d_1+d_2+\alpha r_1)^2}\\
    &=\dfrac{1}{c_\alpha}|Z_\alpha(\cl T)|.
\end{align*}
That is what we want.
\end{proof}
\begin{theorem}\label{finite counting}
    Let $\alpha\in\mathbb{R}$. The stability condition $\sigma_\alpha=(Z_\alpha,\mathcal{P}_\alpha)$ on $\mathcal{T}_C$ satisfies the finite counting condition.
\end{theorem}
\begin{proof}
    We are going to show that for any $\phi\in\mathbb{R}$ and $\gamma\in\mathbb{Z}^4$, the set
    \[\mathcal{P}_\alpha(\phi)_\gamma=\{0\neq T\in\Iso(\mathcal{P}_\alpha(\phi))|\cl T=\gamma\}\]
    is finite. After an integer shift, we may assume that $\phi\in(0,1]$. Let $T=(\mathcal{F}\xrightarrow{f}\mathcal{G})\in\mathcal{P}_\alpha(\phi)_\gamma$, and $\gamma=(r_1,d_1,r_2,d_2)$, where $r_1,r_2\geqslant0$ and when $r_i=0$, $d_i\geqslant0,i=1,2$. Let $\nu=\nu_\alpha(T)$. We only need to show that $\mathcal{F}$ and $\mathcal{G}$ have only finitely many choices up to isomorphisms, and then $f$ also has only finitely many choices up to isomorphisms by the Hom-finiteness of $\Coh C$.

    When $\phi=1$, we must have $r_1=r_2=0,d_1,d_2\geqslant0$. Now each nonzero vertex sheaf is torsion, and hence semistable. By Proposition \ref{semistable sheaf finite} and the fact that all zero sheaf are isomorphic, $\mathcal{F}$ and $\mathcal{G}$ have only finitely many choices up to isomorphisms.

    When $\phi\in(0,1)$, we must have $r_1+r_2>0$. If $r_1=0$, then $\mathcal{F}$ is torsion and hence has only finitely many choices up to isomorphisms. If $r_1>0$, let $\mathcal{E}_1,\dots,\mathcal{E}_n$ be the HN factors of $\mathcal{F}$. By Lemma \ref{slope bounds}(2) and (3), we have
    \[\nu-\alpha\leqslant\mu(\mathcal{E}_i)\leqslant\max(\nu,\nu-\alpha),i=1,\dots,n.\]
    Hence $\mathcal{E}_1,\dots,\mathcal{E}_n$ cannot be torsion, so $\rk\mathcal{E}_i>0,i=1,\dots,n$. Note that $\rk\mathcal{E}_1+\dots+\rk\mathcal{E}_n=r_1$, and there are only finitely many ordered decompositions of $r_1$ into positive integers, so the ranks of HN factors of $\mathcal{F}$ have only finitely many choices. Once the ranks are fixed, the degree of each $\mathcal{E}_i$ also has finitely many choices, since
    \[\nu-\alpha\leqslant\mu(\mathcal{E}_i)=\dfrac{\deg\mathcal{E}_i}{\rk\mathcal{E}_i}\leqslant\max(\nu,\nu-\alpha),\]
    and $\deg\mathcal{E}_i\in\mathbb{Z},i=1,\dots,n$. Thus, by Proposition \ref{semistable sheaf finite}, the HN factors of $\mathcal{F}$ have only finitely many choices up to isomorphisms. Since $\mathcal{F}$ is an iterated extension of its HN factors, and $\Ext^1_{\Coh C}(\mathcal{H},\mathcal{J})$ is finite-dimensional over $k$, and hence finite, for any $\mathcal{H},\mathcal{J}\in\Coh C$, $\mathcal{F}$ has only finitely many choices up to isomorphisms. Similarly, $\mathcal{G}$ also has only finitely many choices up to isomorphisms. Hence we finish the proof.
\end{proof}
\begin{theorem}\label{continuity}
    The map
    \begin{equation}\label{path}
        \begin{split}
            \mathbb{R}&\rightarrow\Stab_\cl(\mathcal{T}_C)\\
            \alpha&\mapsto\sigma_\alpha
        \end{split}
    \end{equation}
    is continuous. As a consequence, every $\sigma_\alpha$ lies in a common connected component $\Sigma_C$ of $\Stab_\cl(\mathcal{T}_C)$.
\end{theorem}
\begin{proof}
    We prove the continuity of (\ref{path}) by showing that $d(\sigma_\alpha,\sigma
    _\beta)\rightarrow0$ when $|\alpha-\beta|\rightarrow0$.

    Let $0\neq T\in\mathcal{P}_\alpha(\phi),\phi\in\mathbb{R}$. Let $S_1,\dots,S_n$ be the HN factors of $T$ with respect to $\sigma_\beta$. Let $m$ be the unique integer such that $T[m]\in\TCoh C$. Then $S_1[m],\dots,S_n[m]$ are the HN factors of $T[m]$ with respect to $\sigma_\beta$. Since $T[m]\in\TCoh C=\mathcal{P}_\beta(0,1]$, we have $S_1[m],\dots,S_n[m]\in\mathcal{P}_\beta(0,1]=\TCoh C$, and $S_1[m]$ is a subobject of $T[m]$, while $S_n[m]$ is a quotient object of $T[m]$. By the $\sigma_\alpha$-semistability of $T[m]$, we have $\theta_\alpha(S_1[m])\leqslant\theta_\alpha(T[m])=\phi+m\leqslant\theta_\alpha(S_n[m])$. By Lemma \ref{phase control}, we obtain
    \begin{align*}
        \phi_\beta^+(T[m])&=\theta_\beta(S_1[m])\leqslant\theta_\alpha(S_1[m])+\dfrac{|\alpha-\beta|}{\pi}\leqslant\phi+m+\dfrac{|\alpha-\beta|}{\pi},\\
        \phi_\beta^-(T[m])&=\theta_\beta(S_n[m])\geqslant\theta_\alpha(S_n[m])-\dfrac{|\alpha-\beta|}{\pi}\geqslant\phi+m-\dfrac{|\alpha-\beta|}{\pi}.
    \end{align*}
    Thus,
    \begin{align}
        \phi_\beta^+(T)&\leqslant\phi+\dfrac{|\alpha-\beta|}{\pi},\label{phi+}\\
        \phi_\beta^-(T)&\geqslant\phi-\dfrac{|\alpha-\beta|}{\pi}.\label{phi-}
    \end{align}
    Now let $0\neq X\in\mathcal{T}_C$. By applying (\ref{phi+}) and (\ref{phi-}) to each $\sigma_\alpha$-HN factor of $X$, we obtain
    \[\phi_\beta^+(X)\leqslant\phi_\alpha^+(X)+\dfrac{|\alpha-\beta|}{\pi},\phi_\beta^-(X)\geqslant\phi_\alpha^-(X)-\dfrac{|\alpha-\beta|}{\pi}.\]
    We can switch $\alpha$ and $\beta$. So
    \begin{equation}\label{phi+-}
        |\phi_\beta^+(X)-\phi_\alpha^+(X)|\leqslant\dfrac{|\alpha-\beta|}{\pi}, |\phi_\beta^-(X)-\phi_\alpha^-(X)|\leqslant\dfrac{|\alpha-\beta|}{\pi}.
    \end{equation}

    Now we focus on the masses of an object with respect to two stability conditions $\sigma_\alpha$ and $\sigma_\beta$. We only consider the case that $|\alpha-\beta|<\dfrac{\pi}{2}$. Let $0\neq T\in\mathcal{P}_\alpha(\phi),\phi\in\mathbb{R}$ with $\sigma_\beta$-HN factors $S_1,\dots,S_n$. By (\ref{phi+}) and (\ref{phi-}), we have
    \[\phi+\dfrac12>\phi+\dfrac{|\alpha-\beta|}{\pi}\geqslant\phi_\beta(S_1)>\dots>\phi_\beta(S_n)\geqslant\phi-\dfrac{|\alpha-\beta|}{\pi}>\phi-\dfrac12.\]
    So
    \[\re(e^{-i\pi\phi}Z_\beta(S_i))=|Z_\beta(S_i)|\cos((\phi_\beta(S_i)-\phi)\pi)\geqslant|Z_\beta(S_i)|\cos|\alpha-\beta|,i=1,\dots,n.\]
    Hence,
    \begin{equation}\label{beta}
        \begin{split}
            |Z_\beta(T)|=\left|\sum_{i=1}^n Z_\beta(S_i)\right|&=\left|\sum_{i=1}^n e^{-i\pi\phi}Z_\beta(S_i)\right|\\
            &\geqslant\re\left(\sum_{i=1}^n e^{-i\pi\phi}Z_\beta(S_i)\right)\\
            &=\sum_{i=1}^n\re(e^{-i\pi\phi}Z_\beta(S_i))\\
            &\geqslant\sum_{i=1}^n|Z_\beta(S_i)|\cos|\alpha-\beta|\\
            &=m_\beta(T)\cos|\alpha-\beta|.
        \end{split}
    \end{equation}
    Note that there exists a unique $m\in\mathbb{Z}$ such that $T[m]\in\mathcal{P}_\alpha(0,1]=\TCoh C$. Write $\cl T[m]=(r_1,d_1,r_2,d_2)$. Then
    \begin{equation}\label{alphabeta}
        \begin{split}
            |Z_\alpha(T)-Z_\beta(T)|&=|Z_\alpha(T[m])-Z_\beta(T[m])|\\
            &=|\alpha-\beta|r_1\\
            &\leqslant|\alpha-\beta||Z_\alpha(T[m])|\\
            &=|\alpha-\beta||Z_\alpha(T)|
        \end{split}
    \end{equation}
    By combining (\ref{beta}) and (\ref{alphabeta}), we obtain
    \begin{equation}\label{semistable mass}
        m_\beta(T)\leqslant\dfrac{1+|\alpha-\beta|}{\cos|\alpha-\beta|}m_\alpha(T).
    \end{equation}
    For any $0\neq X\in\mathcal{T}_C$, by applying (\ref{semistable mass}) for all $\sigma_\alpha$-HN factors $A_1,\dots,A_m$ and summing all the inequalities, we obtain
    \[m_\beta(X)\leqslant\sum_{j=1}^m m_\beta(A_j)\leqslant\dfrac{1+|\alpha-\beta|}{\cos|\alpha-\beta|}\sum_{j=1}^m m_\alpha(A_j)=\dfrac{1+|\alpha-\beta|}{\cos|\alpha-\beta|}m_\alpha(X),\]
    where the first inequality follows from the triangle inequality for mass; see \cite[Proposition 3.3]{Ike2021}, where we take $t=0$. We can switch $\alpha$ and $\beta$, so
    \[m_\alpha(X)\leqslant\dfrac{1+|\alpha-\beta|}{\cos|\alpha-\beta|}m_\beta(X).\]
    Thus,
    \begin{equation}\label{log}
        \left|\log\dfrac{m_\alpha(X)}{m_\beta(X)}\right|\leqslant\log\dfrac{1+|\alpha-\beta|}{\cos|\alpha-\beta|}.
    \end{equation}
    When $|\alpha-\beta|\rightarrow0$, by the definition of the distance and (\ref{phi+-}), (\ref{log}), we have $d(\sigma_\alpha,\sigma_\beta)\rightarrow0$. So the map in (\ref{path}) is continuous.
\end{proof}
\section{Derived Hall algebras and Joyce elements}
\subsection{Derived Hall algebras}
 Let $\mathcal{D}\neq0$ be an essentially small, idempotent-complete $k$-linear triangulated category satisfying the following finiteness conditions as in \cite{XX2008}:
\begin{itemize}
    \item $\mathcal{D}$ is Hom-finite;
    \item $\mathcal{D}$ is Krull-Schmidt, and
    \item $\mathcal{D}$ is left homologically finite, that is, for any $A,B\in\mathcal{D}$
    \[\sum_{i=0}^\infty\dim_k\Hom_{\mathcal{D}}(A[i],B)<+\infty.\]
\end{itemize}
\begin{lemma}
    $\mathcal{T}_C$ satisfies the requirements above.
\end{lemma}
\begin{proof}
    The category $\mathcal{T}_C=D^b(\TCoh C)$ is nonzero, essentially small and $k$-linear. Recall from the proof of Theorem \ref{finite counting} that the $\Hom$ and $\Ext^1$ spaces in $\Coh C$ are finite-dimensional over $k$. For $T,T'\in\TCoh C$, use the short exact sequence (\ref{decomposition}) for both $T$ and $T'$. As in the proof of \cite[Proposition~2.7.24]{Rom2018}, the resulting long exact sequences reduce the finiteness of $\Ext^j_{\TCoh C}(T,T')$ to that of these spaces. Also, $\Ext^j_{\TCoh C}(T,T')=0$ for $j>2$ by Proposition \ref{homdim2}.
    
    Applying the standard truncation triangles to $X,Y\in\mathcal{T}_C$, we obtain
    \[\sum_{i=0}^\infty\dim_k\Hom_{\mathcal{T}_C}(X[i],Y)<+\infty.\]
    Thus $\mathcal{T}_C$ is Hom-finite and left homologically finite. Since idempotents split in the bounded derived category of an abelian category, $\mathcal{T}_C$ is also idempotent-complete. Its Hom-finiteness then gives the Krull--Schmidt property.
\end{proof}

For $A,B\in\mathcal{D}$, put
\[\{A,B\}=\prod_{i>0}|\Hom_{\mathcal{D}}(A[i],B)|^{(-1)^i},a_B=|\Aut_{\mathcal{D}}(B)|.\]

The set of maps $X\rightarrow B$ with cone isomorphic to $A[1]$ is denoted by $\Hom(X,B)_{A[1]}$. Let
\[\mathcal{H}(\mathcal{D})=\bigoplus_{X\in\Iso(\mathcal{D})}\mathbb{Q}u_X,\]
and we define
\begin{equation}\label{eq:Hall}
 u_A*u_B=\sum_X F^X_{AB}u_X,\qquad
 F^X_{AB}=\frac{|\Hom(X,B)_{A[1]}|}{a_B}
              \frac{\{X,B\}}{\{B,B\}}.
\end{equation}
The rational number $F^X_{AB}$ is called the \emph{derived Hall number} with respect to the objects $A,B$ and $X$. 
By \cite[Theorem 3.6]{XX2008}, $(\mathcal{H}(\mathcal{D}),*)$ is an associative $\mathbb{Q}$-algebra with the identity $1=u_0$, called the \emph{derived Hall algebra} of $\mathcal{D}$.
\begin{lemma}[{\cite[Lemma 3.9]{GX2026}}]\label{truncation hall number 1}
    Let $X_{>c}\rightarrow X\rightarrow X_{\leqslant c}\rightarrow X_{>c}[1]$ be the phase truncation of $X$ at $c\in\mathbb{R}$ with respect to some stability condition $\sigma$ on $\mathcal{D}$. Then the derived Hall number
    \[F_{X_{>c}X_{\leqslant c}}^X=1.\]
\end{lemma}
For an element
\[f=\sum_{i=1}^m c_i u_{E_i}\in\mathcal{H}(\mathcal{D}),\]
we can view $f$ as a $\mathbb{Q}$-valued function on $\Iso(\mathcal{D})$:
\begin{align*}
    f:\Iso(\mathcal{D})&\rightarrow\mathbb{Q}\\
    E&\mapsto\begin{cases}
        c_i,& E=E_i,i=1,\dots,m\\
        0,&\text{ otherwise}.
    \end{cases}
\end{align*}
We will frequently use this point of view in this paper. It is clear that any element in $\mathcal{H}(\mathcal{D})$, as a function, is finitely supported.

Let $\mathcal{C}$ be a subcategory of $\mathcal{D}$, with $\Iso(\mathcal{C})$ finite. We denote by $\1_{\mathcal{C}}$ the element in $\mathcal{H}(\mathcal{D})$ taking the value $1$ on objects in $\mathcal{C}$ and $0$ elsewhere. 

$\mathcal{H}(\mathcal{D})$ has a natural grading with respect to $K_0(\mathcal{D})$:
\[\mathcal{H}(\mathcal{D})=\bigoplus_{\gamma\in K_0(\mathcal{D})}\mathcal{H}(\mathcal{D})_\gamma,\]
where $\mathcal{H}(\mathcal{D})_\gamma$ is the subspace of $\mathcal{H}(\mathcal{D})$ spanned by $u_E$ with $E\in\mathcal{D}$ and $\udv E=\gamma$.

Let $\mathcal{C}$ be a subcategory of $\mathcal{D}$. We set the subspace of $\mathcal{H}(\mathcal{D})$:
\[\mathcal{H}(\mathcal{C})=\bigoplus_{X\in\Iso(\mathcal{C})}\mathbb{Q}u_X.\]
When $\mathcal{C}$ is closed under extension, $\mathcal{H}(\mathcal{C})$ is a subalgebra of $\mathcal{H}(\mathcal{D})$. The following lemma provides a useful method for explicit computations:
\begin{lemma}\label{subheart Ringel-Hall}
    Let $\mathcal{B}$ be an abelian subcategory of $\TCoh C$ closed under extension. Then $\mathcal{H}(\mathcal{B})\subseteq\mathcal{H}(\mathcal{T}_C)$ coincides with the Ringel-Hall algebra of $\mathcal{B}$.
\end{lemma}
\begin{proof}
    We refer the reader to \cite[Proposition 6.1]{Toe2006}, but we also give a direct proof.

    Let $U,V\in\mathcal{B}\subseteq\TCoh C$. For $X\in\mathcal{B}$, we set
    \[g_{UV}^X=\#\{U'\subseteq X|U'\cong U,X/U'\cong V\}.\]
    Since $\mathcal{B}$ is closed under extension, we only need to show that
    $F_{UV}^X=g_{UV}^X$
    for any $U,V,X\in\mathcal{B}$. For any $A,B\in\mathcal{B}\subseteq\TCoh C$, we have
    \[\Hom_{\mathcal{T}_C}(A[j],B)=\Hom_{\mathcal{T}_C}(A,B[-j])=0\]
    when $j>0$. So $\{A,B\}=1$ for any $A,B\in\mathcal{B}$. Hence
    \[F_{UV}^X=\dfrac{|\Hom(X,V)_{U[1]}|}{a_V}\dfrac{\{X,V\}}{\{V,V\}}=\dfrac{|\Hom(X,V)_{U[1]}|}{a_V}=\dfrac{\#\{p:X\twoheadrightarrow V|\Ker p\cong U\}}{a_V}.\]
    Now consider the $\Aut V$-action on the set $\{p:X\twoheadrightarrow V|\Ker p\cong U\}$ defined as $\theta\cdot p:=\theta\circ p$. It is free by the surjectivity of $p$. It is easy to see that two maps $p_1,p_2$ in the same orbit must have the same kernel. Conversely, if $\Ker p_1=\Ker p_2=U'$, then $p_1$ and $p_2$ are both cokernels of the inclusion $i:U'\hookrightarrow X$, so there must be a unique $\theta\in\Aut V$ such that $p_2=\theta\circ p_1$, and hence $p_1$ and $p_2$ lie in the same orbit. We also note that for any $U'\subseteq X$ such that $U'\cong U,X/U'\cong V$, the map $p:X\twoheadrightarrow X/U'\xrightarrow{\sim} V$ has kernel $U'\cong U$. So the orbits of $\{p:X\twoheadrightarrow V|\Ker p\cong U\}$ under the $\Aut V$-action are in bijection with the subobjects $U'\subseteq X$ such that $U'\cong U,X/U'\cong V$. By the freeness of the action, we obtain
    \[\#\{p:X\twoheadrightarrow V|\Ker p\cong U\}=a_V\#\{U'\subseteq X|U'\cong U,X/U'\cong V\}=a_Vg_{UV}^X.\]
    Thus we finish the proof.
\end{proof}
\subsection{Completion of fixed slope}
Fix $\mu\in\mathbb{Q}\cup\{\infty\}$. We set
\[\mathcal{B}_\mu=\Mor(\mathcal{C}_\mu).\]
Note that $\mathcal{B}_\mu$ is an abelian subcategory of $\TCoh C$. In fact, we have
\begin{lemma}\label{Bmu}
    $\mathcal{B}_\mu=\mathcal{P}_0(\phi_\mu)$, where $\phi_\mu=\dfrac{1}{\pi}\arg(-\mu+i)\in(0,1)$ when $\mu\in\mathbb{Q}$ and $\phi_\infty=1$.
\end{lemma}
\begin{proof}
    When $\mu\in\mathbb{Q}$, the result follows from Lemma \ref{slope bounds}(4). When $\mu=\infty$, the result follows from the arguments before Lemma \ref{slope bounds}.
\end{proof}

By Lemma \ref{subheart Ringel-Hall}, the subalgebra $\mathcal{H}(\mathcal{B}_\mu)$ coincides with the Ringel-Hall algebra of $\mathcal{B}_\mu$. For $T=(\mathcal{F}\xrightarrow{f}\mathcal{G})\in\TCoh C$, we define the \emph{total rank} of $T$ to be $\rk T=\rk\mathcal{F}+\rk\mathcal{G}$. Let $\mu\in\mathbb{Q}$. For $n\in\mathbb{N}$, consider the subspace
\[\mathcal{H}(\mathcal{B}_\mu)_n=\bigoplus_{T\in\Iso(\mathcal{B}_\mu),\rk T=n}\mathbb{Q}u_T.\]
Then we obtain an $\mathbb{N}$-grading of $\mathcal{H}(\mathcal{B}_\mu)$. For a torsion triple $T=(\mathcal{F}\xrightarrow{f}\mathcal{G})\in\mathcal{B}_\infty$, we define its \emph{total degree} to be $\deg T=\deg\mathcal{F}+\deg\mathcal{G}\geqslant0$. We set
\[\mathcal{H}(\mathcal{B}_\infty)_n=\bigoplus_{T\in\Iso(\mathcal{B}_\infty),\deg T=n}\mathbb{Q}u_T.\]
Then we obtain an $\mathbb{N}$-grading of $\mathcal{H}(\mathcal{B}_\infty)$.
\begin{lemma}
    For any $n\in\mathbb{N}$, the $\mathbb{Q}$-vector space $\mathcal{H}(\mathcal{B}_\mu)_n$ is finite dimensional.
\end{lemma}
\begin{proof}
    We first consider the case $\mu\in\mathbb{Q}$. When $n=0$, note that a triple $T$ with a total rank $0$, with each sheaf having a slope $\mu$, must be zero, so $\mathcal{H}(\mathcal{B}_\mu)_0=\mathbb{Q}u_0$.

    When $n>0$, the number of pairs $(r_1,r_2)$ of natural numbers such that $r_1+r_2=n$ is finite, and for each such pair $(r_1,r_2)$, a triple $T=(\mathcal{F}\xrightarrow{f}\mathcal{G})\in\mathcal{B}_\mu$ with $\rk\mathcal{F}=r_1,\rk\mathcal{G}=r_2$ must satisfy
    \[\deg\mathcal{F}=\mu r_1,\deg\mathcal{G}=\mu r_2,\]
    and $\mathcal{F}$ and $\mathcal{G}$ are both semistable as long as they are not zero. So by Proposition \ref{semistable sheaf finite}, we have only finitely many choices for $\mathcal{F}$ and $\mathcal{G}$, and also for $f\in\Hom_{\mathcal{C}_\mu}(\mathcal{F},\mathcal{G})=\Hom_{\Coh C}(\mathcal{F},\mathcal{G})$ by the Hom-finiteness of $\Coh C$. Thus, the number of isomorphism classes of triples $T\in\mathcal{B}_\mu$ with total rank $n$ is finite, and hence $\mathcal{H}(\mathcal{B}_\mu)_n$ is finite dimensional.

    For the case $\mu=\infty$, one can replace the total ranks with total degrees and apply similar arguments.
\end{proof}

We define the completed algebra
\[\hat{\mathcal{H}}(\mathcal{B}_\mu)=\prod_{n\in\mathbb{N}}\mathcal{H}(\mathcal{B}_\mu)_n.\]
We also define the Lie subalgebra
\[\mathfrak{g}_\mu=\prod_{n\in\mathbb{Z}_+}\mathcal{H}(\mathcal{B}_\mu)_n,\]
with its pro-nilpotent algebraic group $G_\mu=u_0+\mathfrak{g}_\mu$. We define the following two bijections, which are the inverses of each other:
\begin{align*}
    \exp:\mathfrak{g}_\mu&\rightarrow G_\mu&\log:G_\mu&\rightarrow\mathfrak{g}_\mu\\
    x&\mapsto\sum_{i=0}^\infty\dfrac{x^i}{i!}&y&\mapsto\sum_{i=1}^\infty\dfrac{(-1)^{i-1}}{i}(y-u_0)^i.
\end{align*}
We write $\1_{\mathcal{B}_\mu}\in G_\mu\subseteq\hat{\mathcal{H}}(\mathcal{B}_\mu)$ for the characteristic function supported on $\Iso(\mathcal{B}_\mu)$. We now give two decompositions of $\1_{\mathcal{B}_\mu}$, which satisfy a pentagon identity like the one in \cite[\S 4.6]{Bri2018}.
\begin{proposition}\label{ADB}
    Let $0\neq T=(\mathcal{F}\xrightarrow{f}\mathcal{G})\in\mathcal{B}_\mu,\mu\in\mathbb{Q}$. Set
    \[\mathcal{K}=\Ker f,\mathcal{I}=\im f,\mathcal{Q}=\Coker f\in\mathcal{C}_\mu.\]
    When $\alpha\in\mathbb{R}$ changes, the nonzero entries in the following table are the $\sigma_\alpha$-HN factors of $T$, in decreasing order of phase.
    \[\begin{array}{c|c|c}
 \alpha<0&\alpha=0&\alpha>0\\ \hline
 B_{\mathcal{G}},\ A_{\mathcal{F}}&T&A_{\mathcal{K}},\ D_{\mathcal{I}},\ B_{\mathcal{Q}}
 \end{array}.\]
\end{proposition}
\begin{proof}
    The case $\alpha=0$ follows from Lemma \ref{Bmu}. When $\alpha<0$, since $A_\mathcal{F},B_\mathcal{G}$ are $Z_\alpha$-semistable with $Z_\alpha$-slopes $\mu+\alpha$ and $\mu$ respectively when they are nonzero, the short exact sequence (\ref{decomposition}) provides the HN filtration of $T$ and the result follows. When $\alpha>0$, consider
    \begin{equation}\label{+filtration}
        0\subseteq A_\mathcal{K}\subseteq(\mathcal{F}\twoheadrightarrow\mathcal{I})\subseteq T.
    \end{equation}
    The factors of the filtration (\ref{+filtration}) are $A_\mathcal{K},D_\mathcal{I},B_\mathcal{Q}$, with $Z_\alpha$-slope $\mu+\alpha,\mu+\dfrac{\alpha}{2},\mu$ respectively when they are nonzero. It is easy to check that $A_\mathcal{K}$ and $B_\mathcal{Q}$ are $\sigma_\alpha$-semistable when they are nonzero. When $\mathcal{I}\neq0$, we now show that $D_\mathcal{I}$ is $Z_\alpha$-semistable, and hence $\sigma_\alpha$-semistable. Let $E=(\mathcal{U}\xrightarrow{g}\mathcal{V})$ be a nonzero subobject of $D_\mathcal{I}$. Then $g$ is a restriction of $\id_\mathcal{I}$, so it must be injective. Hence, we may assume $\mathcal{U}\subseteq\mathcal{V}\subseteq\mathcal{I}$ and $E=(\mathcal{U}\hookrightarrow\mathcal{V})$. Let $u=\rk\mathcal{U},v=\rk\mathcal{V}$, so $u\leqslant v$. Since $\mathcal{I}$ is semistable, we obtain
    \[\deg\mathcal{U}+\deg\mathcal{V}\leqslant\mu\rk\mathcal{U}+\mu\rk(\mathcal{V})=\mu(u+v).\]
    So by the definition of $Z_\alpha$-slope, we have
    \[\nu_\alpha(E)\leqslant\mu+\alpha\dfrac{u}{u+v}\leqslant\mu+\dfrac{\alpha}{2}=\nu_\alpha(D_\mathcal{I}).\]
    Thus, $D_\mathcal{I}$ is $Z_\alpha$-semistable. Therefore, (\ref{+filtration}) is the $\sigma_\alpha$-HN filtration, and $A_\mathcal{K},D_\mathcal{I},B_\mathcal{Q}$ are the $\sigma_\alpha$-HN factors of $T$ when they are nonzero.
\end{proof}
We define the following three series
\[\Phi_A=\sum_{\mathcal{F}\in\Iso(\mathcal{C}_\mu)}u_{A_\mathcal{F}},\Phi_B=\sum_{\mathcal{F}\in\Iso(\mathcal{C}_\mu)}u_{B_\mathcal{F}},\Phi_D=\sum_{\mathcal{F}\in\Iso(\mathcal{C}_\mu)}u_{D_\mathcal{F}}\in G_\mu\subseteq\hat{\mathcal{H}}(\mathcal{B}_\mu),\]
and
\[L_A=\log\Phi_A,L_B=\log\Phi_B,L_D=\log\Phi_D\in\mathfrak{g}_\mu.\]
\begin{theorem}[Pentagon identity]\label{pentagon}
    We have the identity
    \begin{equation}\label{pentagon identity}
        \1_{\mathcal{B}_\mu}=\Phi_B*\Phi_A=\Phi_A*\Phi_D*\Phi_B
    \end{equation}
    in $G_\mu\subseteq\hat{\mathcal{H}}(\mathcal{B}_\mu)$, and the identity
    \begin{equation}\label{log identity}
        L_D=\log(\exp(-L_A)*\exp(L_B)*\exp(L_A)*\exp(-L_B)).
    \end{equation}
\end{theorem}
\begin{proof}
    For any $0\neq T=(\mathcal{F}\xrightarrow{f}\mathcal{G})\in\mathcal{B}_\mu$, set $\mathcal{K}=\Ker f,\mathcal{I}=\im f,\mathcal{Q}=\Coker f\in\mathcal{C}_\mu$.
    
    Suppose first that $\mu\in\mathbb{Q}$. When $\alpha<0$, taking the phase truncation of $T$ at the phase corresponding to $\mu+\dfrac{\alpha}{2}$, and applying Lemma \ref{truncation hall number 1}, we obtain
    \[(u_{B_\mathcal{G}}*u_{A_\mathcal{F}})(T)=1.\]
    By comparing the coefficient at each  $T\in\mathcal{B}_\mu$, we obtain
    \[\1_{\mathcal{B}_\mu}=\Phi_B*\Phi_A.\]
    
    When $\alpha>0$, taking the phase truncations of $T$ at the phase corresponding to $\mu+\dfrac{\alpha}{4}$ and $\mu+\dfrac{3\alpha}{4}$, and applying Lemma \ref{truncation hall number 1}, we obtain
    \[(u_{A_\mathcal{K}}*u_{D_\mathcal{I}}*u_{B_\mathcal{Q}})(T)=1.\]
    By comparing the coefficient at each $T\in\mathcal{B}_\mu$, we obtain
    \[\1_{\mathcal{B}_\mu}=\Phi_A*\Phi_D*\Phi_B.\]

    It remains to consider the case $\mu=\infty$. In this case every nonzero torsion triple has phase $1$ for all $\alpha\in\mathbb{R}$. We therefore
    compute the two Hall products directly in $\mathcal{B}_\infty$, using Lemma \ref{subheart Ringel-Hall}.
    
    First, $(\Phi_B*\Phi_A)(T)$ counts the subobjects of $T$ of the form $B_\mathcal{V}=(0\rightarrow\mathcal{V})$ whose quotient is of the form $A_\mathcal{U}$. The second vertex of the quotient is $\mathcal{G}/\mathcal{V}$, so we must have $\mathcal{V}=\mathcal{G}$. Thus the only such subobject is $B_\mathcal{G}$, and its quotient is $A_\mathcal{F}$, as in (\ref{decomposition}). Hence $(\Phi_B*\Phi_A)(T)=1$.
    
    Next, $(\Phi_A*\Phi_D*\Phi_B)(T)$ counts the flags
    \[0\subseteq T_1\subseteq T_2\subseteq T\]
    whose successive quotients are of the forms $A$, $D$ and $B$, in this order. Since $T_1$ has zero second vertex and $T/T_2$ has zero first
    vertex, we may write
    \[T_1=(\mathcal{U}\rightarrow0),T_2=(\mathcal{F}\rightarrow\mathcal{V}),\]
    where $\mathcal{U}\subseteq\Ker f$ and $\im f\subseteq\mathcal{V}\subseteq\mathcal{G}$. The arrow of $T_2/T_1$ is the map
    \[\mathcal{F}/\mathcal{U}\longrightarrow\mathcal{V}\]
    induced by $f$. Its kernel is $(\Ker f)/\mathcal{U}$ and its image is $\im f$. Thus $T_2/T_1$ is of the form $D$ if and only if $\mathcal{U}=\Ker f$ and $\mathcal{V}=\im f$. It follows that the only flag counted by the product is
    \[0\subseteq A_\mathcal{K}\subseteq(\mathcal{F}\twoheadrightarrow\mathcal{I})\subseteq T.\]
    Its successive quotients are $A_\mathcal{K}$, $D_\mathcal{I}$ and $B_\mathcal{Q}$. Hence $(\Phi_A*\Phi_D*\Phi_B)(T)=1$. This also covers the case where one of the quotients is zero, since each of the three series has constant term $u_0$. At $T=0$ both products also have coefficient $1$. Comparing coefficients gives
    \[\1_{\mathcal{B}_\infty}=\Phi_B*\Phi_A=\Phi_A*\Phi_D*\Phi_B.\]
    
    Hence, we finish the proof of (\ref{pentagon identity}), so we have
    \[\Phi_D=\Phi_A^{-1}*\Phi_B*\Phi_A*\Phi_B^{-1}=\exp(-L_A)*\exp(L_B)*\exp(L_A)*\exp(-L_B).\]
    Taking the logarithm of both sides, we obtain (\ref{log identity}).
\end{proof}
\subsection{Joyce elements and Joyce transports}
Let $\mathcal{D}$ be a triangulated category satisfying the conditions in Subsection 4.1 with $\cl:K_0(\mathcal{D})\rightarrow\Gamma$ fixed. Let $\sigma=(Z,\mathcal{P})\in\Stab_\cl(\mathcal{D})$ satisfy the support property and the finite counting condition. For any $\phi\in\mathbb{R}$, we define
\[\Lambda_{\sigma,\phi}=\{\cl E|0\neq E\in\mathcal{P}(\phi)\}.\]
For $0\neq\gamma\in\Gamma$, we write
\[\delta^\sigma_{\phi,\gamma}=\sum_{\substack{0\neq E\in\Iso(\mathcal{P}(\phi))\\\cl E=\gamma}}u_E.\]
It is easy to see that $\delta^\sigma_{\phi,\gamma}\neq0$ only if $\gamma\in\Lambda_{\sigma,\phi}$. The sum is finite by the finite counting condition, and so $\delta^\sigma_{\phi,\gamma}\in\mathcal{H}(\mathcal{D})$ for any $\phi\in\mathbb{R}$ and $0\neq\gamma\in\Gamma$.
\begin{definition}[{\cite[\S 4.1]{GX2026}}]
    Given $\phi\in\mathbb{R}$ and $0\neq\gamma\in\Gamma$, define the \emph{Joyce element} to be
    \begin{equation}\label{Joyce element}
        \epsilon^\sigma_{\phi,\gamma}=\sum_{m=1}^\infty\dfrac{(-1)^{m-1}}{m}\sum_{\substack{\gamma_1,\dots,\gamma_m\in\Lambda_{\sigma,\phi}\\\gamma_1+\dots+\gamma_m=\gamma}}\delta^\sigma_{\phi,\gamma_1}*\dots*\delta^\sigma_{\phi,\gamma_m}.
    \end{equation}
\end{definition}
One can easily check that the sum on the right-hand side of (\ref{Joyce element}) must be finite, and hence all Joyce elements lie in $\mathcal{H}(\mathcal{D})$. We denote by $\mathfrak{L}_\sigma$ the Lie subalgebra of $\mathcal{H}(\mathcal{D})$ generated by Joyce elements $\epsilon^\sigma_{\phi,\gamma}$, where $\phi$ runs through $\mathbb{R}$ and $\gamma$ runs through $\Gamma-\{0\}$.

\begin{lemma}[{\cite[Lemma 4.13]{GX2026}}]\label{one point to the component}
    Let $\Sigma\subseteq\Stab_\cl(\mathcal{D})$ be a connected component. If there is a $\sigma\in\Sigma$ that satisfies the support property and the finite counting condition, then any $\tau\in\Sigma$ satisfies the support property and the finite counting condition.
\end{lemma}

Lemma \ref{one point to the component} tells us that, if some $\sigma\in\Sigma$ satisfies the support property and the finite counting condition, then for any $\tau\in\Sigma$, we can define Joyce elements with respect to $\tau$, and then define the Lie subalgebra $\mathfrak{L}_\tau$. Moreover, we have

\begin{theorem}[Joyce transport,{\cite[Theorem 4.14]{GX2026}}]\label{Joyce transport}
    Let $\sigma,\tau\in\Sigma$ where $\Sigma$ is a connected component of $\Stab_\cl(\mathcal{D})$ admitting a stability condition satisfying the support property and the finite counting condition. For any $\phi\in\mathbb{R},0\neq\gamma\in\Gamma$, the Joyce element $\epsilon^\tau_{\phi,\gamma}$ can be written as a finite Lie polynomial in $\epsilon^\sigma_{\psi,\beta},\psi\in\mathbb{R},0\neq\beta\in\Gamma$ with coefficients in $\mathbb{Q}$, with each Lie monomial having total degree $\gamma$. As a consequence, we have
    \[\mathfrak{L}_\sigma=\mathfrak{L}_\tau.\]
\end{theorem}
\begin{corollary}\label{elliptic triple transport}
    The results of Theorem \ref{Joyce transport} hold for $\Sigma_C$. In particular, the results hold for $\sigma_{\alpha_0},\sigma_{\alpha_1}$ for any $\alpha_0,\alpha_1\in\mathbb{R}$. For any $\alpha\in\mathbb{R}$, we have
    \[\mathfrak{L}_{\sigma_\alpha}=\mathfrak{L}_{\sigma_0}.\]
    We will denote the common Lie subalgebra by $\mathfrak{L}_C$.
\end{corollary}
\begin{proof}
    The corollary follows from Theorem \ref{support}, Theorem \ref{finite counting}, Theorem \ref{continuity} and Theorem \ref{Joyce transport}.
\end{proof}
In the rest of this paper, we will focus on the calculation of Joyce elements.
\section{Computations for Joyce elements and the Lie subalgebra}
In this section, we compute Joyce elements of the stability condition $\sigma_0$, and give an explicit form of the Lie subalgebra $\mathfrak{L}_C$. We first introduce a general method using flags, and then apply this method to our computation.
\subsection{Flag polynomials}
Let $\mathcal{B}$ be a $k$-linear Hom-finite essentially small abelian category of finite length. We assume that $\mathcal{B}$ admits an additive $\mathbb{N}$-grading $d$ such that $d(M)>0$ for every $0\neq M\in\mathcal{B}$. We also require that for every $n\in\mathbb{N}$, the number of isomorphism classes of objects in $\mathcal{B}$ of grading $n$ is finite. For example, for any $\mu\in\mathbb{Q}\cup\{\infty\}$, the category $\mathcal{B}_\mu$ defined in Subsection 4.2 satisfies those conditions. We write $\mathcal{H}(\mathcal{B})$ for the Ringel-Hall algebra of $\mathcal{B}$ and $\hat{\mathcal{H}}(\mathcal{B})$ for its completion, constructed via the same process as in Subsection 4.2. We also introduce the exponential and logarithm maps with the same forms.
\begin{lemma}\label{subobject finite}
    For any $0\neq M\in\mathcal{B}$, $M$ admits only finitely many subobjects.
\end{lemma}
\begin{proof}
    For any subobject $U\subseteq M$, we have $d(M)=d(U)+d(M/U)$, so $0\leqslant d(U)\leqslant d(M)$. Since there are only finitely many isomorphism classes in each degree, there are only finitely many choices for $U$ up to isomorphism.

    Fix a representative $U$ of one of these classes. Every subobject of $M$ isomorphic to $U$ is the image of a monomorphism $U\rightarrow M$. Since $\Hom_{\mathcal{B}}(U,M)$ is finite-dimensional over the finite field $k$, there are only finitely many such maps. Thus $M$ has only finitely many subobjects.
\end{proof}
Fix $0\neq M\in\mathcal{B}$. For $m\in\mathbb{N}$, a \emph{weak flag} of $M$ with $m$ steps is a chain $0=F_0\subseteq\cdots\subseteq F_m=M$. We often write it as $F=(F_j)_{j=0}^m$. It is called a \emph{strict flag} if $0=F_0\subsetneqq\cdots\subsetneqq F_m=M$. By Lemma \ref{subobject finite}, for each natural number $m$, the number of $m$-step weak flags of $M$ is finite, and so is that of $m$-step strict flags. We denote by $f_m(M)\in\mathbb{N}$ the number of $m$-step strict flags of $M$. Note that if $m$ is larger than the length of $M$, then $f_m(M)=0$.
\begin{definition}
    We define the \emph{flag polynomial} of $0\neq M\in\mathcal{B}$ to be
    \[Z_M(t)=\sum_{j=1}^{\ell(M)}f_j(M)\binom{t}{j}\in\mathbb{Q}[t],\]
    where $\ell(M)$ is the length of $M$ and
    \[\binom{t}{j}:=\dfrac{t(t-1)\cdots(t-j+1)}{j!}.\]
    When $M=0$, we define $Z_0(t)=1$.
\end{definition}
\begin{lemma}\label{ZMm}
    Let $0\neq M\in\mathcal{B}$ and $m\in\mathbb{N}$. The number of $m$-step weak flags of $M$ is $Z_M(m)$.
\end{lemma}
\begin{proof}
    Let $n\leqslant m$. An $n$-step strict flag gives $\binom mn$ weak flags: choose the positions of its nonzero quotients and add zero quotients on other positions. Conversely, by removing repetitions, one recovers the strict flag. So $Z_M(m)$ counts the $m$-step weak flags.
\end{proof}
One can easily check the following lemma using induction:
\begin{lemma}\label{induction}
    For any $0\neq M\in\mathcal{B}$, we have $f_1(M)=1$ and when $m\geqslant2$,
    \[f_m(M)=\sum_{0\neq U\subsetneqq M}f_{m-1}(U).\]
\end{lemma}
\begin{proposition}\label{derivate at 0}
    For any $0\neq M\in\mathcal{B}$, we have
    \[(\log\1_\mathcal{B})(M)=Z_M'(0)=\sum_{j=1}^{\ell(M)}\dfrac{(-1)^{j-1}}{j}f_j(M).\]
\end{proposition}
\begin{proof}
    Let us first show the second identity. A direct calculation gives
    \begin{align*}
        Z_M'(0)&=\left.\dfrac{\diff}{\diff t}\right|_{t=0}\sum_{j=1}^{\ell(M)}f_j(M)\dfrac{t(t-1)\cdots(t-j+1)}{j!}\\
        &=\sum_{j=1}^{\ell(M)}f_j(M)\left.\dfrac{\diff}{\diff t}\right|_{t=0}\dfrac{t(t-1)\cdots(t-j+1)}{j!}\\
        &=\sum_{j=1}^{\ell(M)}f_j(M)\dfrac{(0-1)\cdots(0-j+1)}{j!}\\
        &=\sum_{j=1}^{\ell(M)}f_j(M)\dfrac{(-1)^{j-1}(j-1)!}{j!}\\
        &=\sum_{j=1}^{\ell(M)}f_j(M)\dfrac{(-1)^{j-1}}{j}.
    \end{align*}
    Now let $x=\1_\mathcal{B}-u_0$. By induction on $j$ and using Lemma \ref{induction}, we obtain
    \[x^j(M)=f_j(M).\]
    Thus, by the definition of $\log$, we have
    \[(\log\1_\mathcal{B})(M)=\sum_{j=1}^\infty\dfrac{(-1)^{j-1}}{j}x^j(M)=\sum_{j=1}^{\ell(M)}\dfrac{(-1)^{j-1}}{j}f_j(M)=Z_M'(0).\]
\end{proof}
\subsection{A formula for \texorpdfstring{$\sigma_0$}{sigma0}-Joyce elements}
We first clarify which Joyce elements we need to calculate. By (\ref{Joyce element}), the Joyce element $\epsilon^{\sigma_0}_{\phi,\gamma}\neq0$ only if the following conditions hold:
\begin{itemize}
    \item $\pi\phi$ is an argument of $Z_0(\gamma)$.
    \item $\gamma$ can be written as the sum of some numerical classes of $\sigma_0$-semistable objects of phase $\phi$. In particular, there exists a $\sigma_0$-semistable object $T$ such that $\gamma=\cl T$.
\end{itemize}
By Lemma \ref{Bmu}, any $\sigma_0$-semistable object lies in $\mathcal{B}_\mu[n]$, where $\mu\in\mathbb{Q}\cup\{\infty\},n\in\mathbb{Z}$. If we set
\[w_\mu=\begin{cases}
    (r,d),&\mu=\dfrac{d}{r},r>0,\gcd(r,d)=1,\\
    (0,1),&\mu=\infty,
\end{cases}\]
then the conditions above can be rewritten as: There exists a $\mu\in\mathbb{Q}\cup\{\infty\}$ such that
$\gamma=(-1)^n(aw_\mu,bw_\mu)$ and $\phi=\phi_\mu+n$ for some $a,b\in\mathbb{N}$ with $a+b>0$ and $n\in\mathbb{Z}$.

Let $F$ be an exact $k$-linear autoequivalence of $\mathcal{T}_C$. We denote by $F_*$ the automorphism of $\mathcal{H}(\mathcal{T}_C)$ induced by $F$, that is, sending each $u_E$ to $u_{F(E)}$. One can easily show that, for any $\gamma=(aw_\mu,bw_\mu)$, we have
\[\epsilon^{\sigma_0}_{\phi_\mu+n,(-1)^n\gamma}=[n]_*\epsilon^{\sigma_0}_{\phi_\mu,\gamma}.\]
Thus, we only need to calculate the Joyce elements of the form
\[\epsilon^{\sigma_0}_{\phi_\mu,(aw_\mu,bw_\mu)},a,b\in\mathbb{N},a+b>0.\]
For $T\in\mathcal{B}_\mu$, if $\cl T=(aw_\mu,bw_\mu)$, then we write $h(T)=a+b$. By Proposition \ref{derivate at 0}, we obtain
\begin{proposition}\label{derivate at 0, Bmu}
    For any $0\neq T\in\mathcal{B}_\mu$ with $\cl T=(aw_\mu,bw_\mu)$, we have
    \[\epsilon^{\sigma_0}_{\phi_\mu,(aw_\mu,bw_\mu)}(T)=(\log\1_{\mathcal{B}_\mu})(T)=Z_T'(0)=\sum_{j=1}^{a+b}\dfrac{(-1)^{j-1}}{j}f_j(T),\]
    where $f_j(T)$ is the number of $j$-step strict flags of $T$.
\end{proposition}
We can write the formula more precisely. Let $\mathcal{S}_1\not\cong\mathcal{S}_2$ be two simple objects in $\mathcal{C}_\mu$, with $T(\mathcal{S}_1),T(\mathcal{S}_2)$ the corresponding stable tubes. Since there is no nontrivial morphism and extension between $T(\mathcal{S}_1)$ and $T(\mathcal{S}_2)$, any morphism and extension between $\Mor(T(\mathcal{S}_1))$ and $\Mor(T(\mathcal{S}_2))$ can only be trivial. If $T=T_1\oplus T_2\in\mathcal{B}_\mu$, where $0\neq T_1\in\Mor(T(\mathcal{S}_1))$ and $0\neq T_2\in\Mor(T(\mathcal{S}_2))$, then for any $m>0$, an $m$-step weak flag of $T$ must be of the form
\[F=(F_{1j}\oplus F_{2j})_{j=0}^m,\]
where $F_i=(F_{ij})_{j=0}^m$ is an $m$-step weak flag of $T_i$ for $i=1,2$, since a subobject of $T$ must be of the form $U_1\oplus U_2$ where $U_1\subseteq T_1,U_2\subseteq T_2$. By Lemma \ref{ZMm}, we have
\[Z_T(t)=Z_{T_1}(t)Z_{T_2}(t).\]
Note that $t|Z_{T_1}(t),t|Z_{T_2}(t)$, so $t^2|Z_T(t)$ and hence $Z_T'(0)=0$, and hence the coefficient of $u_T$ in the Joyce element must be zero. Therefore, we obtain
\begin{theorem}\label{sigma0 formula}
    Let $\mu\in\mathbb{Q}\cup\{\infty\}$ and $a,b\in\mathbb{N}$ with $a+b>0$. Then
    \[\epsilon^{\sigma_0}_{\phi_\mu,(aw_\mu,bw_\mu)}=\sum_{\mathcal{S}\in\Iso(\mathcal{C}_\mu)\text{ simple}}\sum_{\substack{T\in\Iso(\Mor(T(\mathcal{S})))\\
    \cl T=(aw_\mu,bw_\mu)}}Z_T'(0)u_T.\]
\end{theorem}
\subsection{Computation of coefficients}
To compute $Z_T'(0)$ in Theorem \ref{sigma0 formula} is not an easy task. We will only compute it for a triple $T=(\mathcal{F}\xrightarrow{f}\mathcal{G})$ when $\mathcal{F}$ and $\mathcal{G}$ lie in the same stable tube with each nonzero endpoint indecomposable. In this case, by Lemma \ref{Tabl} and the arguments before it, we may assume
\[T=T_{a,b,l}(\mathcal{S})=(J_a(\mathcal{S})\xrightarrow{f_l}J_b(\mathcal{S})),\]
where $\mathcal{S}$ is stable, $a,b\in\mathbb{N},a+b>0$, $0\leqslant l\leqslant\min(a,b)$ and $f_l$ is defined as in (\ref{fl}) when $l\geqslant1$ with $f_0=0$. An easy observation gives
\[T_{a,b,0}(\mathcal{S})=A_{J_a(\mathcal{S})}\oplus B_{J_b(\mathcal{S})}.\]
\begin{lemma}
    When $a,b>0$ and $1\leqslant l\leqslant\min(a,b)$, the triple $T_{a,b,l}(\mathcal{S})$ is indecomposable.
\end{lemma}
\begin{proof}
    Let $T_{a,b,l}(\mathcal{S})=T_1\oplus T_2=(\mathcal{F}_1\xrightarrow{g_1}\mathcal{G}_1)\oplus(\mathcal{F}_2\xrightarrow{g_2}\mathcal{G}_2)$. Then $\mathcal{F}_1\oplus\mathcal{F}_2=J_a(\mathcal{S})$, $\mathcal{G}_1\oplus\mathcal{G}_2=J_b(\mathcal{S})$. By the indecomposability of $J_a(\mathcal{S})$ and $J_b(\mathcal{S})$, we have $\mathcal{F}_1=0$ or $\mathcal{F}_2=0$, and the other is $J_a(\mathcal{S})$. The same holds for $\mathcal{G}_1,\mathcal{G}_2$. If $\mathcal{F}_1=\mathcal{G}_2=0$, or $\mathcal{F}_2=\mathcal{G}_1=0$, we must have $g_1=g_2=0$, and hence $f_l=g_1\oplus g_2=0$, a contradiction. So we must have $\mathcal{F}_1=\mathcal{G}_1=0$ or $\mathcal{F}_2=\mathcal{G}_2=0$, that is, $T_1=0$ or $T_2=0$. So $T_{a,b,l}(\mathcal{S})$ is indecomposable.
\end{proof}
By (\ref{composition series}) and the definition of a subobject in a morphism category, one can easily check that
\begin{lemma}\label{subobjects of Tabl}
    Let $a,b\in\mathbb{N},a+b>0$ and $0\leqslant l\leqslant\min(a,b)$. A subobject of $T_{a,b,l}(\mathcal{S})$ must be of the form
    \[T_{u,v,l(u)}(\mathcal{S}),0\leqslant v\leqslant b,0\leqslant u\leqslant\min(a,v+a-l),\]
    where $l(u)=\max(0,u-a+l)$. Conversely, a triple of the above form is indeed a subobject of $T_{a,b,l}(\mathcal{S})$. Thus, the subobjects of $T_{a,b,l}(\mathcal{S})$ are in bijection with pairs $(u,v)$ such that $0\leqslant v\leqslant b,0\leqslant u\leqslant\min(a,v+a-l)$, and $T_{u_1,v_1,l(u_1)}(\mathcal{S})\subseteq T_{u_2,v_2,l(u_2)}(\mathcal{S})$ if and only if $u_1\leqslant u_2$ and $v_1\leqslant v_2$.
\end{lemma}
We set
\[W_{a,b,l}(t)=Z_{T_{a,b,l}(\mathcal{S})}(t),N=a+b-l.\]
Then we have
\begin{lemma}\label{Wabl}
    Let $m\in\mathbb{Z}_+$. When $l>0$, we have
    \[W_{a,b,l}(m)=\binom{m+a-1}{a}\binom{m+b-1}{b}-\binom{m+N-1}{N+1}\binom{m+l-1}{l-1},\]
    and when $l=0$, we have
    \[W_{a,b,0}(m)=\binom{m+a-1}{a}\binom{m+b-1}{b}.\]
\end{lemma}
\begin{proof}
    By Lemma \ref{ZMm}, $W_{a,b,l}(m)$ is equal to the number of $m$-step weak flags of $T_{a,b,l}(\mathcal{S})$. By Lemma \ref{subobjects of Tabl}, it is equal to the number of sequences $(u_0,v_0),\dots,(u_m,v_m)$, such that
    \begin{equation}\label{condition 1}
        0=u_0\leqslant\dots\leqslant u_m=a,0=v_0\leqslant\dots\leqslant v_m=b,
    \end{equation}
    and
    \begin{equation}\label{condition 2}
        u_j\leqslant v_j+a-l,j=1,\dots,m.
    \end{equation}
    The number of such sequences satisfying (\ref{condition 1}) is
    \[\binom{m+a-1}{a}\binom{m+b-1}{b}.\]
    
    When $l=0$, the condition (\ref{condition 2}) holds naturally, since $u_j\leqslant a\leqslant v_j+a$. So we have
    \[W_{a,b,0}(m)=\binom{m+a-1}{a}\binom{m+b-1}{b}.\]
    
    When $l>0$, we have to count the sequences that violate (\ref{condition 2}). We set $s=a-l$ for convenience. For each sequence $(u_0,v_0),\dots,(u_m,v_m)$ satisfying (\ref{condition 1}), we encode it by two northeast lattice paths in $\mathbb{Z}^2$, that is, paths going only to the right or upwards: the first path $P$ (the blue one in Figure \ref{fig:weak-flag-lattice-paths}) is from $A_1=(0,0)$ to $B_1=(a,m-1)$, with horizontal run $u_j-u_{j-1}$ at height $j-1$, while the second path $R$ (the orange one in Figure \ref{fig:weak-flag-lattice-paths}) is from $A_2=(s+1,-1)$ to $B_2=(s+1+b,m-2)$, with horizontal run $v_j-v_{j-1}$ at height $j-2$.
    
    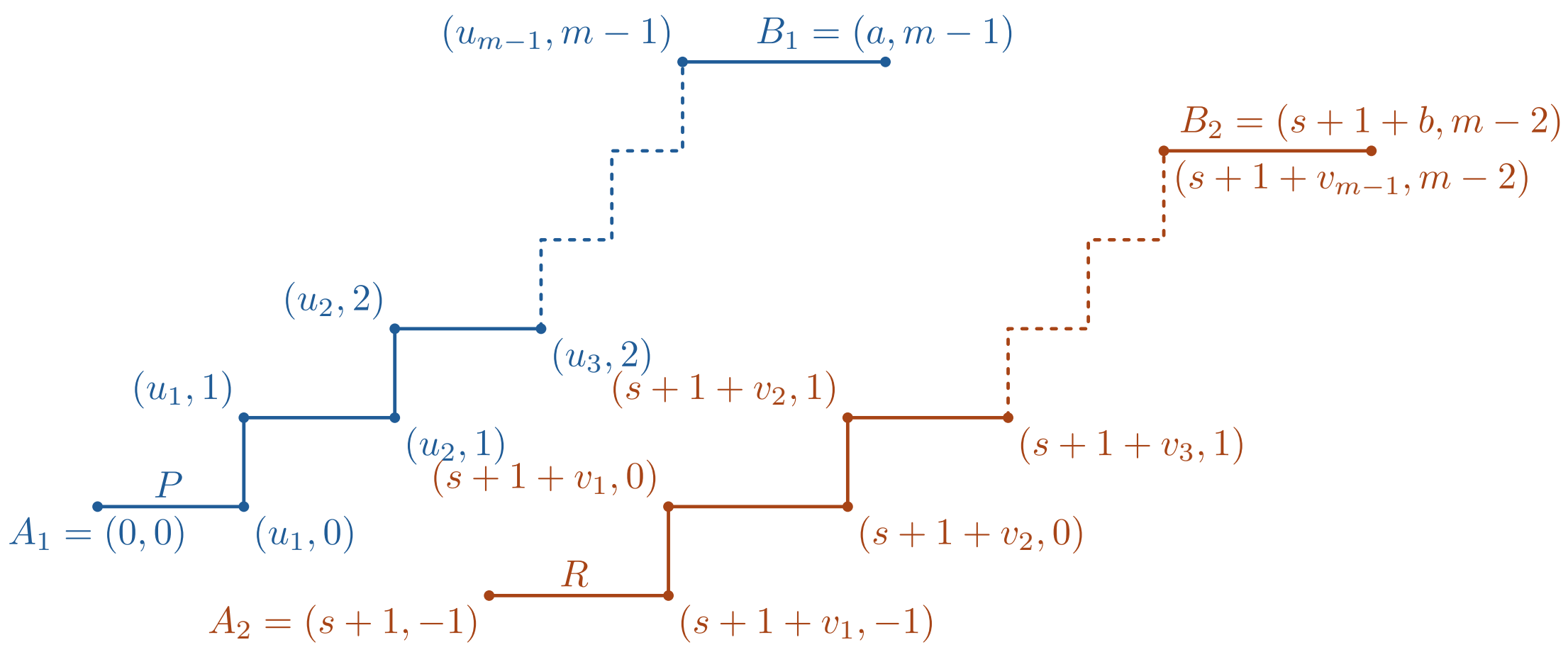
\begin{figure}[htbp]
\centering
\begin{tikzpicture}[
    x=.88cm, y=.83cm, font=\footnotesize,
    line width=.9pt, line cap=round, line join=round,
    every node/.style={inner sep=2pt},
    omitted/.style={dash pattern=on 1.3pt off 2.6pt}
]
  \definecolor{PathBlue}{HTML}{205C97}
  \definecolor{PathOrange}{HTML}{A74416}
  \def\pathmark#1#2#3{
    \fill #1 circle[radius=1.4pt];
    \node[#2] at #1 {$#3$};}
  \begin{scope}[color=PathBlue]
    \draw (0,0) -- (1.55,0) -- (1.55,1) -- (3.15,1)
                -- (3.15,2) -- (4.70,2);
    \draw[omitted]
          (4.70,2) |- (5.45,3) |- (6.20,4) -- (6.20,5);
    \draw (6.20,5) -- (8.35,5);

    \node[above] at (.75,0) {$P$};
    \pathmark{(0,0)}   {below}       {A_1=(0,0)}
    \pathmark{(1.55,0)}{below right} {(u_1,0)}
    \pathmark{(1.55,1)}{above left}  {(u_1,1)}
    \pathmark{(3.15,1)}{below right} {(u_2,1)}
    \pathmark{(3.15,2)}{above left}  {(u_2,2)}
    \pathmark{(4.70,2)}{below right} {(u_3,2)}
    \pathmark{(6.20,5)}{above left}  {(u_{m-1},m-1)}
    \pathmark{(8.35,5)}{above}       {B_1=(a,m-1)}
  \end{scope}
  \begin{scope}[color=PathOrange,shift={(4.15,-1)}]
    \draw (0,0) -- (1.90,0) -- (1.90,1) -- (3.80,1)
                -- (3.80,2) -- (5.50,2);
    \draw[omitted]
          (5.50,2) |- (6.35,3) |- (7.15,4) -- (7.15,5);
    \draw (7.15,5) -- (9.35,5);

    \node[above] at (.90,0) {$R$};
    \pathmark{(0,0)}   {below left}  {A_2=(s+1,-1)}
    \pathmark{(1.90,0)}{below right} {(s+1+v_1,-1)}
    \pathmark{(1.90,1)}{above left}  {(s+1+v_1,0)}
    \pathmark{(3.80,1)}{below right} {(s+1+v_2,0)}
    \pathmark{(3.80,2)}{above left}  {(s+1+v_2,1)}
    \pathmark{(5.50,2)}{below right} {(s+1+v_3,1)}
    \pathmark{(7.15,5)}{below right} {(s+1+v_{m-1},m-2)}
    \pathmark{(9.35,5)}{above}       {B_2=(s+1+b,m-2)}
  \end{scope}
\end{tikzpicture}
\caption{The lattice paths associated with a weak flag.}
\label{fig:weak-flag-lattice-paths}
\end{figure}
    
    The sequences satisfying (\ref{condition 1}) are in one-to-one correspondence with the pairs $(P,R)$ of northeast paths $P:A_1\rightarrow B_1,R:A_2\rightarrow B_2$. The condition (\ref{condition 2}) is satisfied if and only if the two paths do not meet. Thus our task is to count the pairs $(P,R)$ such that $P$ and $R$ meet. We claim that the number of such pairs is equal to the number of pairs $(P',R')$ of northeast paths $P':A_1\rightarrow B_2,R':A_2\rightarrow B_1$. In fact, if $P$ and $R$ meet, let $H$ be their first common point, and we construct two northeast paths:
    \[P':A_1\xrightarrow{P}H\xrightarrow{R}B_2,R':A_2\xrightarrow{R}H\xrightarrow{P}B_1.\]
    Conversely, note that two northeast paths $P':A_1\rightarrow B_2$ and $R':A_2\rightarrow B_1$ must meet. Let $H$ be their first common point. We construct
    \[P:A_1\xrightarrow{P'}H\xrightarrow{R'}B_1,R:A_2\xrightarrow{R'}H\xrightarrow{P'}B_2.\]
    One can check that the two constructions are inverses to each other, so the claim holds. Now the number of northeast paths from $A_1$ to $B_2$ is $\binom{m+N-1}{N+1}$, and the number of northeast paths from $A_2$ to $B_1$ is $\binom{m+l-1}{l-1}$. Thus we obtain
    \[W_{a,b,l}(m)=\binom{m+a-1}{a}\binom{m+b-1}{b}-\binom{m+N-1}{N+1}\binom{m+l-1}{l-1}.\]
\end{proof}
\begin{theorem}\label{epsilon Tabl}
    Let $\mathcal{S}\in\mathcal{C}_\mu$ be stable. Write $e=[K_\mathcal{S}:k]$. Let $a,b\in\mathbb{N}, a+b>0$ and $0\leqslant l\leqslant\min(a,b)$. When $a>0,b>0$, we have
    \[\epsilon^{\sigma_0}_{\phi_\mu,(aew_\mu,bew_\mu)}(T_{a,b,l}(\mathcal{S}))=\begin{cases}
        \dfrac{1}{(a+b-l)(a+b-l+1)},&l>0,\\
        0,&l=0.
    \end{cases}\]
    When $a=0$ or $b=0$, we have
    \[\epsilon^{\sigma_0}_{\phi_\mu,(aew_\mu,bew_\mu)}(T_{a,b,0}(\mathcal{S}))=\dfrac{1}{a+b}.\]
\end{theorem}
\begin{proof}
    By Proposition \ref{derivate at 0, Bmu}, we have
    \[\epsilon^{\sigma_0}_{\phi_\mu,(aew_\mu,bew_\mu)}(T_{a,b,l}(\mathcal{S}))=Z_{T_{a,b,l}(\mathcal{S})}'(0)=W_{a,b,l}'(0).\]
    
    By Lemma \ref{Wabl}, when $l=0$, 
    \[W_{a,b,0}(t)=\binom{t+a-1}{a}\binom{t+b-1}{b}.\]
    When $a>0,b>0$, both factors vanish at $t=0$, so $t$ divides both factors, and hence $t^2$ divides $W_{a,b,0}(t)$, which leads to $W_{a,b,0}'(0)=0$. When $a=0$, the first factor is $1$ and hence $W'_{0,b,0}(0)$ is equal to the derivative of the second factor at $t=0$, which is $\dfrac{1}{b}$. Similarly when $b=0$ we have $W'_{a,0,0}(0)=\dfrac{1}{a}$.

    Also by Lemma \ref{Wabl}, when $l>0$,
    \[W_{a,b,l}(t)=\binom{t+a-1}{a}\binom{t+b-1}{b}-\binom{t+N-1}{N+1}\binom{t+l-1}{l-1}.\]
    By the same argument, the derivative of the first term at $0$ is $0$, so we only need to calculate the derivative of the second term. We also note that the first factor of the second term takes the value $0$ and the second factor takes $1$ when $t=0$, so
    \begin{align*}
        W_{a,b,l}'(0)&=-\left.\dfrac{\diff}{\diff t}\right|_{t=0}\binom{t+N-1}{N+1}\binom{t+l-1}{l-1}\\
        &=-\left.\dfrac{\diff}{\diff t}\right|_{t=0}\binom{t+N-1}{N+1}\\
        &=-\left.\dfrac{\diff}{\diff t}\right|_{t=0}\dfrac{(t+N-1)(t+N-2)\cdots t(t-1)}{(N+1)!}\\
        &=-\dfrac{(0+N-1)(0+N-2)\cdots(0+1)(0-1)}{(N+1)!}\\
        &=\dfrac{(N-1)!}{(N+1)!}=\dfrac{1}{N(N+1)}=\dfrac{1}{(a+b-l)(a+b-l+1)}.
    \end{align*}
\end{proof}
\subsection{The Lie subalgebra}
Fix $\mu\in\mathbb{Q}\cup\{\infty\}$. Let $\mathcal{H}(\mathcal{C}_\mu)$ be the Ringel-Hall algebra of the abelian category $\mathcal{C}_\mu$. We write $\hat{\mathcal{H}}(\mathcal{C}_\mu)$ for its completion, by the same construction as in Subsection 4.2. For a letter $X\in\{A,B,D\}$, consider the linear map defined below and extended continuously:
\begin{align*}
    X_*:\hat{\mathcal{H}}(\mathcal{C}_\mu)&\rightarrow\hat{\mathcal{H}}(\mathcal{B}_\mu)\\
    u_\mathcal{F}&\mapsto u_{X_\mathcal{F}}.
\end{align*}
One can show that $A_*$ and $B_*$ also preserve the Hall multiplication, and hence are algebra homomorphisms.

We write $\epsilon^C_{\mu,(r,d)}$ for the $(r,d)$-component of $\log\1_{\mathcal{C}_\mu}\in\hat{\mathcal{H}}(\mathcal{C}_\mu)$. In particular, we set
\[e_{\mu,m}=\epsilon^C_{\mu,mw_\mu}.\]
Thus we have
\begin{equation}\label{log1=emum}
    \log\1_{\mathcal{C}_\mu}=\sum_{m=1}^\infty e_{\mu,m}.
\end{equation}
\begin{lemma}\label{A*B*}
    For any $m\in\mathbb{Z}_+$ and $\mu\in\mathbb{Q}\cup\{\infty\}$, we have
    \[\epsilon^{\sigma_0}_{\phi_\mu,(mw_\mu,0)}=A_*e_{\mu,m},\epsilon^{\sigma_0}_{\phi_\mu,(0,mw_\mu)}=B_*e_{\mu,m}.\]
\end{lemma}
\begin{proof}
    We prove the first identity. By Lemma \ref{Bmu} and (\ref{Joyce element}), $\epsilon^{\sigma_0}_{\phi_\mu,(mw_\mu,0)}$ is the $(mw_\mu,0)$-component of $\log\1_{\mathcal{B}_\mu}$. Let $T=(\mathcal{F}\xrightarrow{f}\mathcal{G})\in\mathcal{B}_\mu$ be such that $\cl T=(mw_\mu,0)$. Then $\rk\mathcal{G}=\deg\mathcal{G}=0$, so $\mathcal{G}=0$, and hence $T=A_{\mathcal{F}}$ with $\cl\mathcal{F}=mw_\mu$.

    A subobject of $A_{\mathcal{F}}$ in $\mathcal{B}_\mu$ has the form $A_{\mathcal{U}}$, where $\mathcal{U}$ is a subobject of $\mathcal{F}$ in $\mathcal{C}_\mu$. Conversely, each such $\mathcal{U}$ gives a subobject $A_{\mathcal{U}}$ of $A_\mathcal{F}$, with quotient $A_\mathcal{F}/A_\mathcal{U}\cong A_{\mathcal{F}/\mathcal{U}}$. Thus, for any $j\in\mathbb{Z}_+$, the $j$-step strict flags of $A_{\mathcal{F}}$ in $\mathcal{B}_\mu$ are in bijection with those of $\mathcal{F}$ in $\mathcal{C}_\mu$. So the flag polynomial of $\mathcal{F}\in\mathcal{C}_\mu$ is equal to the flag polynomial of $A_\mathcal{F}\in\mathcal{B}_\mu$. By Proposition \ref{derivate at 0} and Proposition \ref{derivate at 0, Bmu}, we obtain
    \[\epsilon^{\sigma_0}_{\phi_\mu,(mw_\mu,0)}(T)=(\log\1_{\mathcal{B}_\mu})(A_{\mathcal{F}})=(\log\1_{\mathcal{C}_\mu})(\mathcal{F})=e_{\mu,m}(\mathcal{F}).\]
    Since $A_*u_{\mathcal{F}}=u_{A_{\mathcal{F}}}$, comparison of the coefficients gives $\epsilon^{\sigma_0}_{\phi_\mu,(mw_\mu,0)}=A_*e_{\mu,m}$. For the class $(0,mw_\mu)$, the first vertex is zero, and the subobjects of $B_{\mathcal{G}}$ are precisely the triples $B_{\mathcal{V}}$ with $\mathcal{V}\subseteq\mathcal{G}$ in $\mathcal{C}_\mu$. The same comparison gives the second identity.
\end{proof}
\begin{lemma}\label{LALB}
    For any $\mu\in\mathbb{Q}\cup\{\infty\}$, we have
    \[L_A=\sum_{m=1}^\infty A_*e_{\mu,m},L_B=\sum_{m=1}^\infty B_*e_{\mu,m}.\]
\end{lemma}
\begin{proof}
    We only prove the first identity, since the proof of the second is similar. Recall that
    \[L_A=\log\Phi_A=\log\left(\sum_{\mathcal{F}\in\Iso(\mathcal{C}_\mu)}u_{A_\mathcal{F}}\right).\]
    Since $A_*:\hat{\mathcal{H}}(\mathcal{C}_\mu)\rightarrow\hat{\mathcal{H}}(\mathcal{B}_\mu)$ is an algebra homomorphism, we have
    \begin{align*}
        L_A&=\log\left(\sum_{\mathcal{F}\in\Iso(\mathcal{C}_\mu)}u_{A_\mathcal{F}}\right)\\
        &=\log\left(\sum_{\mathcal{F}\in\Iso(\mathcal{C}_\mu)}A_*u_{\mathcal{F}}\right)\\
        &=\log\left(A_*\sum_{\mathcal{F}\in\Iso(\mathcal{C}_\mu)}u_{\mathcal{F}}\right)\\
        &=\log(A_*\1_{\mathcal{C}_\mu})\\
        &=A_*(\log\1_{\mathcal{C}_\mu})\\
        (\ref{log1=emum})&=A_*\left(\sum_{m=1}^\infty e_{\mu,m}\right)\\
        &=\sum_{m=1}^\infty A_*e_{\mu,m}.
    \end{align*}
\end{proof}
\begin{theorem}\label{generators of LC}
    The Lie subalgebra $\mathfrak{L}_C$ of $\mathcal{H}(\mathcal{T}_C)$ in Corollary \ref{elliptic triple transport} is equal to
    \[\mathfrak{L}=\Lie_{\mathbb{Q}}\linp{[n]_*A_*e_{\mu,m},[n]_*B_*e_{\mu,m}\middle|\mu\in\mathbb{Q}\cup\{\infty\},m\in\mathbb{Z}_+,n\in\mathbb{Z}}\subseteq\mathcal{H}(\mathcal{T}_C).\]
\end{theorem}
\begin{proof}
    By Corollary \ref{elliptic triple transport}, Lemma \ref{Bmu} and the arguments in Subsection 5.2, we have
    \[\mathfrak{L}_C=\mathfrak{L}_{\sigma_0}=\Lie_{\mathbb{Q}}\linp{[n]_*\epsilon^{\sigma_0}_{\phi_\mu,(aw_\mu,bw_\mu)}\middle|\mu\in\mathbb{Q}\cup\{\infty\},a,b\in\mathbb{N},a+b>0,n\in\mathbb{Z}}.\]
    So $\mathfrak{L}\subseteq\mathfrak{L}_C$ follows immediately from Lemma \ref{A*B*}.

    To show that $\mathfrak{L}_C\subseteq\mathfrak{L}$, we only need to show that for any $a,b\in\mathbb{N},a+b>0$ and $\mu\in\mathbb{Q}\cup\{\infty\}$, the Joyce element $\epsilon^{\sigma_0}_{\phi_\mu,(aw_\mu,bw_\mu)}$ can be written as a finite Lie polynomial in some $A_*e_{\mu,m}$ and $B_*e_{\mu,l}$. Consider the series
    \[\log\1_{\mathcal{B}_\mu}=\sum_{\substack{a,b\in\mathbb{N}\\a+b>0}}\epsilon^{\sigma_0}_{\phi_\mu,(aw_\mu,bw_\mu)}.\]
    By the pentagon identity (\ref{pentagon identity}), we have
    \[\log\1_{\mathcal{B}_\mu}=\log(\Phi_B*\Phi_A)=\BCH(L_B,L_A),\]
    where $\BCH(x,y):=\log(\exp x*\exp y)$ is the Baker-Campbell-Hausdorff Lie series. For any $a,b\in\mathbb{N},a+b>0$, the $(a,b)$-term $\epsilon^{\sigma_0}_{\phi_\mu,(aw_\mu,bw_\mu)}$ of the series $\log\1_{\mathcal{B}_\mu}$ is a Lie polynomial of terms of degrees no more than $a+b$ in $L_B$ and $L_A$. By Lemma \ref{LALB}, the terms of degrees no more than $a+b$ in $L_A$ and $L_B$ are finite, and they are $A_*e_{\mu,m}$ and $B_*e_{\mu,l}$ for some $l,m\leqslant a+b$. Thus the result follows.
\end{proof}
\section*{Acknowledgments}
This work was supported by NSF of China (No. 12371036).
\section*{Declaration of AI use}
The author used ChatGPT to assist in improving some of the mathematical arguments and the English expression. ChatGPT also helped the author draw Figure \ref{fig:weak-flag-lattice-paths}. All mathematical results, arguments, and proofs were carefully developed and verified by the author.

\end{document}